\documentclass[11pt]{article}
\usepackage{amssymb,amsmath,amsthm}
\newcommand{\DD}{\mathcal{D}}
\newcommand{\EE}{\mathcal{E}}

\newcommand{\OO}{\mathcal{O}}

\newcommand{\C}{\mathbb{C}}
\newcommand{\N}{\mathbb{N}}

\newcommand{\R}{\mathbb{R}}
\newcommand{\T}{\mathbb{T}}
\newcommand{\Z}{\mathbb{Z}}

\newcommand{\dd}{\,{\rm d}}
\newcommand{\D}{{\rm d}}

\renewcommand{\div}{\mathop{\mathrm{div}}}
\newcommand{\curl}{\mathop{\mathrm{curl}}}

\renewcommand{\:}{\thinspace :}

\newcommand{\ul}{\mathrm{ul}}

\newcommand{\BUC}{\mathrm{BUC}}

\newtheorem{theorem}{Theorem}[section]

\newtheorem{proposition}[theorem]{Proposition}
\newtheorem{lemma}[theorem]{Lemma}
\newtheorem{corollary}[theorem]{Corollary}

\theoremstyle{definition}
\newtheorem{remark}[theorem]{Remark}
\newtheorem{remarks}[theorem]{Remarks}

\newcommand{\QED}{\mbox{}\hfill$\Box$}
\numberwithin{equation}{section}

\begin{document}
 
\title{Uniform boundedness and long-time asymptotics
for the two-dimensional Navier-Stokes equations in an 
infinite cylinder}

\author{
\null\\
{\bf Thierry Gallay}\\ 
Institut Fourier\\
Universit\'e de Grenoble 1\\
BP 74\\
38402 Saint-Martin-d'H\`eres, France\\
{\tt Thierry.Gallay@ujf-grenoble.fr}
\and
\\
{\bf Sini\v{s}a Slijep\v{c}evi\'c}\\    
Department of Mathematics\\ 
University of Zagreb\\
Bijeni\v{c}ka 30\\
10000 Zagreb, Croatia\\
{\tt slijepce@math.hr}}

\date{}

\maketitle

\begin{abstract}\noindent
We study the incompressible Navier-Stokes equations in the 
two-dimensional strip $\R \times [0,L]$, with periodic 
boundary conditions and no exterior forcing. If the initial
velocity is bounded, we prove that the solution remains uniformly 
bounded for all times, and that the vorticity distribution
converges to zero as $t \to \infty$. We deduce that, after a
transient period, a laminar regime emerges in which the 
solution rapidly converges to a shear flow governed by 
the one-dimensional heat equation. Our approach is constructive
and gives explicit estimates on the size of the solution and 
the lifetime of the turbulent period in terms of the initial 
Reynolds number. 
\end{abstract}

\maketitle

\section{Introduction}\label{secIntro}

We are interested in understanding the dynamics of the incompressible
Navier-Stokes equations in large or unbounded spatial domains. In
particular, for initial data with bounded energy density, we would
like to estimate the kinetic energy of the solution in a small
subdomain at a given time, independently of the total initial energy
which may be infinite if the domain is unbounded. In other words, we
are looking for {\em uniformly local} energy estimates that would
control how much energy can be transfered from one region to another 
in the system.

This question is already interesting in the relatively simple
situation where the fluid is supposed to evolve in a bounded
two-dimensional domain $\Omega \subset \R^2$, with no-slip boundary
condition and no exterior forcing. In that case, if the initial data
are bounded, it is well known that the solution of the Navier-Stokes
equations is globally defined in the energy space and converges to
zero, at an exponential rate, as $t \to +\infty$ \cite{CF}. This
certainly implies that the fluid velocity $u(x,t)$ is uniformly
bounded for all times, but all estimates we are aware of depend on the
size of the domain $\Omega$ or on the total initial energy, and not
only on the initial energy {\em density}. Indeed, although the 
total energy of the fluid is a decreasing function of time, 
the fluid velocity $u(x,t)$ may temporarily increase in some 
regions due to energy redistribution in the system. 

To control these fluctuations, it is rather natural to begin with the
idealized situation where the Navier-Stokes equations are considered
in the whole space $\R^2$, with initial data that are merely
bounded. In that case, it is possible to prove the existence of a
unique global solution \cite{GMS}, provided the pressure is defined in
an appropriate way \cite{Ka}. The corresponding velocity $u(x,t)$
belongs to $L^\infty(\R^2)$ for all $t \ge 0$, but it is not known
whether the norm $\|u(\cdot,t)\|_{L^\infty}$ stays uniformly bounded
for all times. Early results gave pessimistic estimates on that
quantity \cite{GMS,ST}, but a substantial progress was recently made
by S.~Zelik \cite{Ze2} who showed that $\|u(\cdot,t)\|_{L^\infty}$
cannot grow faster than $t^2$ as $t \to \infty$, see
Section~\ref{secConc} for a more precise statement. Still, we do not
have any example of unbounded solution, and it is therefore unclear
whether the above result is optimal.

The aim of the present paper is to address these issues in the 
the simplified setting where the fluid velocity and the pressure
are supposed to be {\em periodic} in one space direction. In other
words, we consider the incompressible Navier-Stokes equations in the 
two-dimensional strip $\Omega_L = \R \times [0,L]$, with periodic 
boundary conditions. Our system reads
\begin{equation}\label{NSeq}
  \partial_t u + (u\cdot \nabla)u \,=\, \nu\Delta u -\frac{1}{\rho}
  \nabla p~, \qquad \div u \,=\, 0~,
\end{equation}
where $u(x,t) \in \R^2$ is the velocity field and $p(x,t) \in \R$ the
associated pressure. We denote the space variable by $x = (x_1,x_2)$,
where $x_1 \in \R$ will be called the ``horizontal'' coordinate and
$x_2 \in [0,L]$ the ``vertical'' coordinate. The physical parameters
in \eqref{NSeq} are the kinematic viscosity $\nu > 0$ and the fluid
density $\rho > 0$, which are both supposed to be constant.  Besides
the pressure, an important quantity derived from the velocity $u$ is
the vorticity $\omega = \partial_1 u_2 - \partial_2 u_1$, which
satisfies the advection-diffusion equation
\begin{equation}\label{Veq}
  \partial_t \omega + u\cdot \nabla \omega \,=\, \nu\Delta \omega~.
\end{equation}

As was already explained, we want to consider infinite-energy
solutions of \eqref{NSeq}, for which the velocity field is merely
bounded.  We thus assume that, for any $t \ge 0$, the velocity
$u(\cdot,t)$ belongs to $\BUC(\Omega_L)$, the space of all bounded and
uniformly continuous functions $u : \Omega_L \to \R^2$ that satisfy
the periodicity condition $u(x_1,0) = u(x_1,L)$ for all $x_1 \in
\R$. It is clear that $\BUC(\Omega_L)$ is a Banach space when equipped
with the uniform norm
\[
  \|u\|_{L^\infty(\Omega_L)} \,=\, \sup_{x \in \Omega_L} |u(x)|~, 
  \qquad \hbox{where} \quad |u| = (u_1^2 + u_2^2)^{1/2}~.
\]
If $u \in \BUC(\Omega_L)$ is divergence-free and if the associated
vorticity $\omega$ is bounded, one can show that the elliptic equation
$-\Delta p = \rho\div(u\cdot \nabla)u$ has a bounded solution 
$p : \Omega_L\to \R$ such that $p(x_1,0) = p(x_1,L)$ for all 
$x_1 \in \R$. Moreover, there exists $C > 0$ such that
\[
  \|p\|_{L^\infty(\Omega_L)} \,\le\, C \rho L^2 
  \|\omega\|_{L^\infty(\Omega_L)}^2~,
\]
see \cite[Lemma~2.3]{GS3} and Section~\ref{subsecP4} below. 
This is our definition of the pressure in \eqref{NSeq}, which 
agrees (up to an irrelevant additive constant) with the 
choice made in \cite{GMS,Ka} in a more general situation.  

Given divergence-free initial data $u_0 \in \BUC(\Omega_L)$ with
associated vorticity distribution $\omega_0$, we introduce the 
Reynolds numbers
\begin{equation}\label{Reyndef}
  R_u \,=\, \frac{L}{\nu}\,\|u_0\|_{L^\infty(\Omega_L)}~, \qquad 
  R_\omega \,=\, \frac{L^2}{\nu}\,\|\omega_0\|_{L^\infty(\Omega_L)}~. 
\end{equation}
The following result shows that the Cauchy problem for \eqref{NSeq} 
is globally well-posed in the space $\BUC(\Omega_L)$. 

\begin{theorem}\label{Cauchy}
For any $u_0 \in \BUC(\Omega_L)$ with $\div u_0 = 0$, the Navier-Stokes 
equations \eqref{NSeq} with the above choice of the pressure have 
a unique global solution $u \in C^0([0,+\infty),\BUC(\Omega_L))$ 
such that $u(0) = u_0$. Moreover, there exists a constant $C > 0$, 
depending only on the initial Reynolds number $R_u$, such that 
\begin{equation}\label{ubound}
  \frac{L}{\nu}\|u(\cdot,t)\|_{L^\infty(\Omega_L)} \,\le\, C
  \Bigl(1 + \frac{\sqrt{\nu t}}{L}\Bigr)~, \qquad 
  \hbox{for all }~ t \ge 0~.
\end{equation}
\end{theorem}

Existence of a unique global solution to \eqref{NSeq} in $\BUC(\Omega_L)$ 
is ensured by the general results of Giga, Matsui, and Sawada
\cite{GMS}, which apply to the Navier-Stokes equations in the whole
plane $\R^2$ with initial data in $\BUC(\R^2)$. The specific situation
where the flow is periodic in one space direction was considered by
Afendikov and Mielke \cite{AM}, with the motivation of understanding
the transition to instability in Kolmogorov flows. In the particular
case where no exterior forcing is applied, the results of \cite{AM}
give an upper bound on $\|u(\cdot,t)\|_{L^\infty}$ which grows
linearly in time, and can be improved with little extra effort to 
provide estimate \eqref{ubound}, see \cite{GS3}. A further progress
was made in \cite{GS3}, where the authors proved that 
$\|u(\cdot,t)\|_{L^\infty}$ cannot grow faster than $t^{1/6}$ as 
$t \to \infty$. Moreover, several results were obtained in 
\cite[Theorem~1.3]{GS3} which strongly suggest that solutions of 
\eqref{NSeq} in $\Omega_L$ should stay uniformly bounded for all times. 
For instance, for all $T > 0$, we have the following estimate
\[
  \sup_{x_1 \in \R}\int_0^T \!\!\int_0^L |u(x_1,x_2,t)|^2 \dd x_2 \dd t 
  \,\le\, C\,\frac{\nu^2 T}{L}~,
\]
for some $C > 0$ depending only on the initial Reynolds number 
$R_u$. In addition, for all $R > 0$ and all $T > 0$, one finds
\begin{equation}\label{uintbound}
  \int_{B_R} |u(x,T)|^2 \dd x + \nu \int_0^T \!\!\int_{B_R} 
  |\nabla u(x,t)|^2 \dd x \dd t \,\le\, C\,\frac{\nu^2}{L}
  \Bigl(R + \sqrt{\nu T}\Bigr)~,
\end{equation}
where $B_R = \{(x_1,x_2) \in \Omega_L\,|\, |x_1| \le R\}$. 
Estimate \eqref{uintbound} shows in particular that the energy 
dissipation rate $\nu |\nabla u(x,t)|^2$ converges to zero on 
average as $t \to \infty$, and this in turn implies that the 
solution $u(x,t)$ approaches for ``almost all'' times the family 
of spatially homogeneous steady states of \eqref{NSeq}, uniformly 
on compact subsets of $\Omega_L$, see \cite[Section~8]{GS3} 
for a precise statement. 

In this paper, we complement and substantially improve the 
results of \cite{GS3} by showing that any solution of 
\eqref{NSeq} in $\BUC(\Omega_L)$ remains {\em uniformly bounded}
for all times, and converges as $t \to \infty$ {\em uniformly on} 
$\Omega_L$ to a simple shear flow of the form
\begin{equation}\label{shear}
  u_\infty(x,t) \,=\, \begin{pmatrix}c \\ m(x_1,t)\end{pmatrix}~,
  \qquad p(x,t) \,=\, 0~,
\end{equation}
where $c \in \R$ is a constant and $m(x_1,t)$ is an approximate 
solution of the one-dimensional heat equation $\partial_t m \,=\, 
\nu\partial_1^2 m$. Our main result can be stated as follows\:

\begin{theorem}\label{main}
For any divergence-free initial data $u_0 \in \BUC(\Omega_L)$ 
with bounded vorticity distribution $\omega_0$, the solution
of the Navier-Stokes equations \eqref{NSeq} given by 
Theorem~\ref{Cauchy} has the following properties\: 
\\[1mm]
{\bf 1.} {\em (Uniform boundedness of the velocity)} There exists $C > 0$ 
such that, for all $t \ge 0$, 
\begin{equation}\label{main1}
  \frac{L}{\nu}\,\|u(t)\|_{L^\infty(\Omega_L)} \,\le\, C\Bigl(R_u + R_\omega 
  + (1+R_\omega)(R_u^2 + R_\omega^2)\Bigr)~,
\end{equation}
where $R_u,R_\omega$ are given by \eqref{Reyndef}.\\[1mm]
{\bf 2.} {\em (Uniform decay of the vorticity)} There exists $C > 0$ 
such that, for all $t > 0$, 
\begin{equation}\label{main2}
  \Bigl(\frac{L^2}{\nu}\,\|\omega(t)\|_{L^\infty(\Omega_L)}\Bigr)^2 
  \,\le\, C(1+R_\omega)(R_u^2 + R_\omega^2)\frac{L}{\sqrt{\nu t}}~.
\end{equation}
{\bf 3.} {\em (Exponential convergence to a shear flow)} For any
$\gamma < 2\pi^2$, we have
\begin{equation}\label{main3}
  \frac{L}{\nu}\,\|u(t)- u_\infty(t)\|_{L^\infty(\Omega_L)} \,=\, 
  \OO\Bigl(\exp{\Bigl(-\frac{\gamma\nu t}{L^2}
  \Bigr)}\Bigr)~, \qquad t \to \infty~.
\end{equation}
Here $u_\infty(x,t)$ is given by \eqref{shear}, where $c \in \R$ 
and $m(x_1,t)$ is an approximate solution of the one-dimensional heat 
equation, in the sense that $\partial_t m - \nu \partial_1^2 m = 
\OO(e^{-2\gamma\nu t/L^2})$ as $t \to \infty$. 
\end{theorem}

\begin{remarks}\label{aftermain}\quad\\
{\bf 1.} The constant $C$ in estimates \eqref{main1}, 
\eqref{main2} is universal\: the dependence of both members upon 
the initial data $u_0$ and the physical parameters $\nu,\rho,L$ 
is entirely explicit. Note also that inequalities 
\eqref{main1}--\eqref{main3} involve only dimensionless quantities, 
such as the initial Reynolds numbers $R_u$ and $R_\omega$. 
\\[1mm]
{\bf 2.} The assumption that the initial vorticity $\omega_0 = 
\curl u_0$ be bounded is by no mean essential. Indeed, due to 
parabolic regularization, any solution of \eqref{NSeq} given by 
Theorem~\ref{Cauchy} is smooth for positive times, and 
has a bounded vorticity distribution for $t > 0$. More quantitatively,
the existence proof given in \cite{GS3} shows that there exists a
positive constant $C$ such that $\|u(\cdot,t_0)\|_{L^\infty} \le 
2 \|u_0\|_{L^\infty}$ and $\nu\|\nabla u(\cdot,t_0)\|_{L^\infty} \le 
C \|u_0\|_{L^\infty}^2$ if $t_0 = C^{-1}\nu\|u_0\|_{L^\infty}^{-2}$. 
The Reynolds numbers \eqref{Reyndef} computed at time $t_0$ 
thus satisfy $R_u(t_0) \le 2R_u$ and $R_\omega(t_0) \le C R_u^2$. 
\\[1mm]
{\bf 3.} It follows from \eqref{main3} that all steady states of the
Navier-Stokes equations in $\BUC(\Omega_L)$ are spatially homogeneous\:
$u = (c_1,c_2)^\top$, with $c_1,c_2 \in \R$. For the same reason,
nontrivial time-periodic solutions or more general recurrent orbits
do not exist. Since bounded solutions of the heat equation on $\R$
converge uniformly on compact sets toward the family of constant 
solutions, we also deduce from \eqref{main3} that all solutions 
of \eqref{NSeq} given by Theorem~\ref{Cauchy} converge uniformly
on compact sets to the family of spatially homogeneous steady 
states as $t \to \infty$. Note that this conclusion is stronger
than what one typically expects for general extended dissipative 
systems, see \cite{GS2}.
\\[1mm]
{\bf 4.} By the parabolic maximum principle, the vorticity bound
$\|\omega(\cdot,t)\|_{L^\infty}$ is nonincreasing with time, and
\eqref{main2} shows that this quantity converges to zero as $t \to
\infty$. As we shall see in Section~\ref{secExp}, when the vorticity
Reynolds number $L^2 \|\omega(\cdot,t)\|_{L^\infty}/\nu$ becomes smaller
than a universal constant (related to the Poincar\'e inequality), 
the system enters a laminar regime where the solution rapidly 
converges to a shear flow. Thus, for large initial data, we can
identify two different stages in the evolution of the system\: 
a long transient period, in which turbulence can develop, and 
an asymptotic laminar regime described by \eqref{main3}. 
In view of \eqref{main2}, the lifetime $T$ of the turbulent 
period satisfies $\nu T/L^2 \le CR^6$ for some $C > 0$, where 
$R = \max(R_u,R_\omega)$. 
\\[1mm]
{\bf 5.} Although our motivation for using periodic boundary 
conditions is to shed some light on the behavior of solutions 
to the Navier-Stokes equations in the whole plane $\R^2$, it 
is natural at this point to ask what happens if we consider 
\eqref{NSeq} in the strip $\Omega_L = \R \times [0,L]$ with 
other conditions at the boundary. If we assume that the velocity 
$u$ vanishes on $\partial\Omega$ ({\em no-slip boundary 
conditions}), then the solution of \eqref{NSeq} decays exponentially 
to zero as $t \to \infty$, see \cite{Ze1,AZ}. We thus have the 
analog of \eqref{main3} with $u_\infty = 0$. Another interesting 
possibility is to suppose that $u_2 = \partial_2 u_1 = 0$ on 
$\partial\Omega$ ({\em perfect slip boundary conditions}). In 
that case, if we extend the solution $u(x,t)$ to the larger strip 
$\R \times [-L,L]$ in such a way that $u_1$ (resp. $u_2$)
is an even (resp. odd) function of the vertical coordinate $x_2$, 
it is easy to verify that the extended velocity field satisfies 
periodic boundary conditions on $\R \times [-L,L]$. In addition, 
the vertical velocity is by construction an odd function of $x_2$, 
hence has a zero vertical average. It follows that \eqref{main1}, 
\eqref{main2} hold, as well as \eqref{main3} with $u_\infty = 
(c,0)^\top$ and $\gamma < \pi^2/2$. Finally, it is also possible to 
consider {\em Navier friction conditions}, but this intermediate 
case has not been studied so far, and our approach does not apply 
directly due to the lack of a priori estimate on the vorticity. 
\end{remarks}

The rest of this paper is devoted to the proof of Theorem~\ref{main},
which relies on previous results from \cite{GS3} and is also strongly
inspired by our recent work on extended dissipative systems
\cite{GS2}.  In Section~\ref{secPrelim} below, we recall a few basic
facts about equation \eqref{NSeq} which were already established in
\cite{AM,GS3}. In particular, we single out the important role played
by the vertical average of the velocity field, which cannot be simply
estimated using the Biot-Savart law and the a priori bound on the
vorticity. We also give an explicit formula for the pressure in
\eqref{NSeq}. In Section~\ref{secVort}, we study in some detail the
linear advection-diffusion equation \eqref{Veq} in $\Omega_L$, with
periodic boundary conditions, assuming that the velocity field
$u(x,t)$ is given. Using ideas that date back from the pioneering work
of Nash \cite{Na,FS}, we establish an accurate upper bound on the
fundamental solution of \eqref{Veq} which shows that, if $\div u = 0$
and if the first component $u_1$ has zero vertical average, solutions
of \eqref{Veq} spread diffusively as $t \to \infty$. The core of the
proof begins in Section~\ref{secEner}, where we control the evolution
of the velocity and vorticity fields using weighted energy
estimates. This is strongly related to the approach developed in
\cite{GS3}, although the new formulation we propose here is 
completely self-contained. Combining the weighted energy and 
enstrophy estimates of Section~\ref{secEner} with the results 
of Section~\ref{secVort}, we prove in Section~\ref{secUnif} the 
first two assertions of Theorem~\ref{main}, namely the uniform 
bound on the velocity field and the decay estimate for the 
vorticity. In Section~\ref{secExp}, we study the time evolution
of small solutions of \eqref{NSeq} and show that they converge 
to shear flows as $t \to \infty$. Finally, some conclusions 
and perspectives are presented in Section~\ref{secConc}, while
Section~\ref{secApp} is an appendix which contains the proof of 
a technical lemma stated in Section~\ref{secVort}. 

\medskip\noindent{\bf Acknowledgements.} The authors thank 
A. Mielke and S. Zelik for fruitful discussions. The research of 
Th.G. is partially support by ANR grant DYFICOLTI from the French 
Ministry of Research. 

\section{Preliminaries}\label{secPrelim}

In this section we recall some basic properties of equation
\eqref{NSeq} which were already established in \cite{AM,GS3}, and we
prepare the proof of Theorem~\ref{main} by performing a 
few preliminary steps.

\subsection{Nondimensionalization}\label{subsecP1}
Our system contains three physical parameters: the kinematic viscosity
$\nu$, the fluid density $\rho$, and the width $L$ of the spatial
domain. All of them can be eliminated if we introduce the new
variables $\tilde x = x/L \in \R \times [0,1]$, $\tilde t =\nu t/L^2
\ge 0$, and the new functions $\tilde u, \tilde \omega,\tilde p$ 
defined by the relations
\begin{equation}\label{adim}
  u(x,t) \,=\, \frac{\nu}{L}\,\tilde u\Bigl(\frac{x}{L},
  \frac{\nu t}{L^2}\Bigr)~, \qquad
  \omega(x,t) \,=\, \frac{\nu}{L^2}\,\tilde \omega\Bigl(\frac{x}{L},
  \frac{\nu t}{L^2}\Bigr)~, \qquad
  p(x,t) \,=\, \frac{\rho\nu^2}{L^2}\,\tilde p\Bigl(\frac{x}{L},
  \frac{\nu t}{L^2}\Bigr)~.
\end{equation}
In what follows, we work exclusively with the rescaled variables 
$\tilde x, \tilde t$ and the dimensionless functions $\tilde u, 
\tilde \omega, \tilde p$, but we drop the tildes for notational 
simplicity. We thus consider the nondimensionalized Navier-Stokes 
equations
\begin{equation}\label{NSeq2}
  \partial_t u + (u\cdot \nabla)u \,=\, \Delta u -\nabla p~, 
  \qquad \div u \,=\, 0~, 
\end{equation}
as well as the corresponding vorticity equation
\begin{equation}\label{Veq2}
  \partial_t \omega + u\cdot \nabla \omega \,=\, \Delta \omega~.
\end{equation}
In both systems, since we impose periodic boundary conditions, it 
is mathematically convenient to assume that the space variable 
$x = (x_1,x_2)$ lies in the cylinder $\Omega = \R \times \T$, where 
$\T = \R/\Z$ is the one-dimensional torus. Using \eqref{adim}, 
it is straightforward to translate Theorems~\ref{Cauchy} and 
\ref{main} into their equivalent, nondimensional form. In particular, 
we observe that the Reynolds numbers defined in \eqref{Reyndef} 
are now simply given by $R_u = \|u_0\|_{L^\infty(\Omega)}$ and 
$R_\omega = \|\omega_0\|_{L^\infty(\Omega)}$. 

\subsection{Decomposition of the velocity}\label{subsecP2}
Let $u(x,t)$ be a solution of the Navier-Stokes equation \eqref{NSeq2} 
given by Theorem~\ref{Cauchy}. Due to the incompressibility condition, 
the vertical average of the horizontal velocity
\begin{equation}\label{average}
  \langle u_1\rangle(x_1,t) \,:=\, \int_\T u_1(x_1,x_2,t)\dd x_2 
\end{equation}
satisfies $\partial_1 \langle u_1\rangle = 0$, and using \eqref{NSeq2}
together with our definition of the pressure one can also show that
$\partial_t \langle u_1\rangle = 0$, see \cite{AM,GS3}. Thus $\langle
u_1\rangle$ is a constant which can be eliminated using an appropriate
Galilean transformation, without affecting our results in any
essential way.  In what follows, we assume therefore that $\langle
u_1\rangle = 0$, so that the velocity $u(x,t)$ has the following
decomposition\:
\begin{equation}\label{udecomp}
  u(x,t) \,=\, \begin{pmatrix} 0 \\ m(x_1,t)\end{pmatrix}
  + \widehat u(x,t)~, \qquad \hbox{where}\quad 
  m(x_1,t) \,=\, \int_\T u_2(x_1,x_2,t)\dd x_2~.
\end{equation}
As is easily verified, the mean vertical speed $m$ and the oscillating 
part $\widehat u = (\widehat u_1,\widehat u_2)^\top$ satisfy the 
evolution equations
\begin{align}\label{meq}
  \partial_t m + \partial_1 \langle \widehat u_1 \widehat u_2\rangle
  \,&=\, \partial_1^2 m~, \\ \label{hatu1eq}
  \partial_t \widehat u_1 + \widehat u_1\partial_1 \widehat u_1 
  + (m+\widehat u_2)\partial_2\widehat u_1 \,&=\, \Delta 
  \widehat u_1 - \partial_1 p~, \\ \label{hatu2eq}
  \partial_t \widehat u_2 + \widehat u_1\partial_1 \widehat u_2 
  + (m+\widehat u_2)\partial_2\widehat u_2 +  \widehat u_1\partial_1 m
  - \partial_1 \langle \widehat u_1 \widehat u_2\rangle
  \,&=\, \Delta \widehat u_2 - \partial_2 p~,
\end{align}
where the brackets $\langle \cdot\rangle$ denote the vertical 
average, as in \eqref{average}. In a similar way, we can decompose
the vorticity as $\omega(x,t) = \partial_1 m(x_1,t) + \widehat 
\omega(x,t)$, where $\langle \widehat\omega \rangle = 0$.  

\subsection{The Biot-Savart law}\label{subsecP3}
As is explained in \cite{AM,GS3}, the oscillating part of the velocity 
can be reconstructed from the vorticity via the Biot-Savart formula 
$\widehat u = \nabla^\perp K*\widehat \omega$, where $\nabla^\perp = 
(-\partial_2,\partial_1)^\top$ and $K$ is the fundamental solution of 
the Laplace operator in $\Omega = \R \times \T$\: 
\begin{equation}\label{Kdef}
  K(x_1,x_2) \,=\, \frac{1}{4\pi}\log\Bigl(2\cosh(2\pi x_1) 
  - 2\cos(2\pi x_2)\Bigr)~, \qquad x \in \R^2~.
\end{equation}
Explicitly, we have
\begin{equation}\label{BS}
  \widehat u(x) \,=\, \int_\Omega\nabla^\perp K(x-y) 
  \widehat\omega(y)\dd y~, \qquad x \in \Omega~.
\end{equation}
In contrast, the mean vertical flow $m(x_1,t)$ cannot be fully
reconstructed from the vorticity, and we only know that $\partial_1 m 
= \langle \omega\rangle$. It follows in particular from \eqref{BS} that
\begin{equation}\label{BSres}
  \|\widehat u\|_{L^\infty(\Omega)} \,\le\, C_1 \|\widehat \omega
  \|_{L^\infty(\Omega)} \,\le\, 2C_1 \|\omega\|_{L^\infty(\Omega)}~,
\end{equation}
for some $C_1 > 0$.

\subsection{Definition of the pressure}\label{subsecP4}
The pressure satisfies the following elliptic equation in $\Omega = 
\R \times \T$\:
\begin{equation}\label{peq2}
  -\Delta p \,=\, \div(u\cdot \nabla)u \,=\, 
  \Delta(u_1^2) + 2 \partial_2 (\omega u_1)~,
\end{equation}
where the second equality is an identity that holds for any  
divergence-free vector field $u$ with $\omega = \partial_1 u_2 - 
\partial_2 u_1$. Since we assume that $p$ is bounded (i.e., there 
is no pressure gradient at infinity), then $p$ is determined by
\eqref{peq2} up to an irrelevant additive constant, and we have
\begin{equation}\label{pdef}
  p \,=\, -u_1^2 - 2\partial_2 K*(\omega u_1)~.
\end{equation}
As $u_1 = \widehat u_1$, it follows from \eqref{BSres} and
\eqref{pdef} that
\begin{equation}\label{pres}
  \|p\|_{L^\infty} \,\le\, C_2 \|\omega\|_{L^\infty}^2~, 
\end{equation}
for some $C_2 > 0$. 

\section{Estimates for the vorticity equation}\label{secVort}
\label{subEst}

In this section, we assume that we are given a smooth divergence-free 
vector field $u(x,t)$ on the two-dimensional cylinder $\Omega = \R 
\times \T$, and we study the linear advection-diffusion equation
\begin{equation}\label{Veq3}
  \partial_t \omega + u\cdot \nabla \omega \,=\, \Delta \omega~, 
  \qquad x \in \Omega~, \quad t \ge 0~.
\end{equation}
Our goal is to establish a priori estimates on the solutions of
\eqref{Veq3}, in the spirit of the fundamental work of Nash
\cite{Na}. These estimates are well known when Eq.~\eqref{Veq3} is
considered in the whole space $\R^n$, but the case of the product
manifold $\Omega = \R \times \T$ is apparently less documented in the
literature (see however \cite{CGL,Gr1,Gr2}). In any case, the proofs are
rather standard, and we reproduce them below for the reader's
convenience.

\subsection{The Nash inequality in $\R \times \T$}

In the whole space $\R^n$, it was shown by Nash \cite{Na}
that there exists a constant $C_n > 0$ such that
\begin{equation}\label{NashRn}
  \|f\|_{L^2(\R^n)}^{1+2/n} \,\le\, C_n \|\nabla f\|_{L^2(\R^n)}
  \|f\|_{L^1(\R^n)}^{2/n}~,
\end{equation}
for any $f \in H^1(\R^n) \cap L^1(\R^n)$. In the cylinder 
$\Omega = \R \times \T$ inequality \eqref{NashRn} does 
not hold, but we have the following estimate, which can be
interpreted as a combination of \eqref{NashRn} with $n = 1$ 
and $n = 2$. 

\begin{lemma}\label{Nash}
There exists a constant $C > 0$ such that, for all $f \in H^1(\Omega) 
\cap L^1(\Omega)$, 
\begin{equation}\label{Nashbd}
  \|f\|_{L^2(\Omega)} \,\le\, C\max\Bigl\{\|\nabla f\|_{L^2(\Omega)}^{1/3}
  \|f\|_{L^1(\Omega)}^{2/3}\,,\, \|\nabla f\|_{L^2(\Omega)}^{1/2}
  \|f\|_{L^1(\Omega)}^{1/2}\Bigr\}~.
\end{equation}
\end{lemma}

\noindent{\bf Proof.} We mimick the proof of the classical Nash 
inequality \cite{Na}. Given a nonzero $f \in H^1(\Omega) \cap 
L^1(\Omega)$, we use the Fourier representation
\[
  f(x) \,=\, \int_{\hat \Omega} \hat f(\xi)\,e^{i\xi\cdot x}\dd \mu(\xi)~, 
  \qquad 
  \hat f(\xi) \,=\, \int_\Omega f(x)\,e^{-i\xi\cdot x}\dd x~, 
\]
where $\hat \Omega = \R \times (2\pi\Z)$ is the dual manifold
and $\mu$ is the positive measure on $\hat \Omega$ defined by 
\[
  \int_{\hat \Omega} \phi(\xi)\dd\mu(\xi) \,=\, \frac{1}{2\pi}
  \int_\R \sum_{n \in \Z} \phi(k,2\pi n)\dd k~,
\]
for any continuous function $\phi : \hat \Omega \to \C$ with 
compact support. Given $R > 0$, the Parseval identity implies
\begin{equation}\label{Nash1}
  \int_\Omega |f(x)|^2 \dd x \,=\, \int_{\hat \Omega} |\hat f(\xi)|^2 
  \dd\mu(\xi) \,=\, \int_{\hat \Omega_R} |\hat f(\xi)|^2 \dd\mu(\xi) 
  + \int_{\hat \Omega_R^c} |\hat f(\xi)|^2 \dd\mu(\xi)~,
\end{equation}
where $\hat\Omega_R \,=\, \{\xi \in \hat \Omega\,|\, |\xi| \le R\}$. 
To estimate the right-hand side of \eqref{Nash1}, we observe that
\begin{align*}
  \int_{\hat \Omega_R} |\hat f(\xi)|^2 \dd\mu(\xi) \,&\le\, 
  \|\hat f\|_{L^\infty(\hat \Omega)}^2\,\mu(\hat \Omega_R) 
  \,\le\, \|f\|_{L^1(\Omega)}^2 \,\mu(\hat \Omega_R)~, \\
  \int_{\hat \Omega_R^c} |\hat f(\xi)|^2 \dd\mu(\xi) \,&\le\, 
  \int_{\hat \Omega_R^c} \frac{|\xi|^2}{R^2}|\hat f(\xi)|^2 \dd\mu(\xi) 
  \,\le\, \frac{1}{R^2} \|\nabla f\|_{L^2(\Omega)}^2~,
\end{align*}
and it is easy to verify that $\mu(\hat \Omega_R) \le \max(R,R^2)$ for 
any $R > 0$. We thus have
\begin{equation}\label{Nash2}
  \|f\|_{L^2(\Omega)}^2 \,\le\, \|f\|_{L^1(\Omega)}^2\,\max(R,R^2)
  + \frac{1}{R^2}\|\nabla f\|_{L^2(\Omega)}^2~.
\end{equation}
If we now choose
\[
  R \,=\, \begin{cases} \|\nabla f\|_{L^2(\Omega)}^{2/3} 
  \|f\|_{L^1(\Omega)}^{-2/3} & \hbox{if}~\, \|\nabla f\|_{L^2(\Omega)} 
  \le \|f\|_{L^1(\Omega)}~, \\[1mm]
  \|\nabla f\|_{L^2(\Omega)}^{1/2} \|f\|_{L^1(\Omega)}^{-1/2}
  & \hbox{if}~\, \|\nabla f\|_{L^2(\Omega)} \ge \|f\|_{L^1(\Omega)}~,
  \end{cases}
\]
we see that \eqref{Nashbd} follows from \eqref{Nash2}. \QED

\bigskip 
We use below an equivalent form of \eqref{Nashbd}, which is
called a $\psi$-Nash inequality in \cite{CGL}\:

\begin{corollary}\label{psiNash}
There exists a constant $C > 0$ such that, for all nonzero 
$f \in H^1(\Omega) \cap L^1(\Omega)$, 
\begin{equation}\label{psiNashbd}
  \|\nabla f\|_{L^2(\Omega)} \,\ge\, C \|f\|_{L^2(\Omega)}
  \min\left\{\frac{\|f\|_{L^2(\Omega)}}{\|f\|_{L^1(\Omega)}}\,,\, 
  \frac{\|f\|_{L^2(\Omega)}^2}{\|f\|_{L^1(\Omega)}^2}\right\}~.
\end{equation}
\end{corollary}

\subsection{$L^p-L^q$ estimates}\label{subLpLq}

Using Nash's inequality, we next derive $L^p-L^q$ estimates 
for solutions of \eqref{Veq3}. 

\begin{proposition}\label{LpLqProp}
Given $1 \le p \le q \le \infty$, there exists a constant $K_1 > 0$ 
(independent of $u$) such that any solution of \eqref{Veq3} with 
initial data $\omega_0 \in L^p(\Omega)$ satisfies
\begin{equation}\label{LpLq}
   \|\omega(t)\|_{L^q(\Omega)} \,\le\, \frac{K_1}{V(t)^{\frac1p-
  \frac1q}} \,\|\omega_0\|_{L^p(\Omega)}~, \qquad t > 0~, 
\end{equation}
where $V(t) = \min(t,\sqrt{t})$. 
\end{proposition}

\begin{remark}\label{volume}
We observe that $V(t)$ is, up to inessential constants, the volume 
of a ball of radius $\sqrt{t}$ in the manifold $\Omega = \R \times 
\T$. The fact that estimate \eqref{LpLq} holds with a constant $K_1$ 
independent of $u$ is intuitively clear, since the advection term 
$u \cdot \nabla \omega$ in \eqref{Veq3} does not affect $L^p$ norms.
\end{remark}

\noindent{\bf Proof.} Since the velocity field $u(x,t)$ in
\eqref{Veq3} is divergence-free, it is well-known that, for any $p \in
[1,\infty]$, the $L^p$ norm of any solution of \eqref{Veq3} is a
nonincreasing function of time. This means that \eqref{LpLq} holds
with $K_1 = 1$ and $q = p$ for any $p \in [1,\infty]$. Thus it remains 
to prove \eqref{LpLq} for $p = 1$, $q = \infty$, and the other cases will
follow by interpolation, see \cite[Section 1.1]{BL}.

We argue as in \cite[Proposition II.1]{Co}. Let $\omega(x,t)$ be 
a solution of \eqref{Veq3} with nonzero initial data $\omega_0 \in 
L^1(\Omega)$. Since the velocity field $u(x,t)$ is divergence-free, 
a direct calculation shows that
\[
 \frac{\D}{\D t}\int_\Omega \omega(x,t)^2\dd x \,=\, -2 \int_\Omega
 |\nabla \omega(x,t)|^2 \dd x \,\le\, 0~, \qquad t > 0~. 
\]
To estimate the right-hand side, we use Nash's inequality 
\eqref{psiNashbd} which gives
\[
  \|\nabla\omega(t)\|_{L^2(\Omega)}^2 \,\ge\, C \|\omega(t)\|_{L^2(\Omega)}^2
  \min\left\{\frac{\|\omega(t)\|_{L^2(\Omega)}^2}{\|\omega_0\|_{L^1(\Omega)}^2}~,\,
  \frac{\|\omega(t)\|_{L^2(\Omega)}^4}{\|\omega_0\|_{L^1(\Omega)}^4}\right\}~,
  \qquad t > 0~,
\]
since $\|\omega(t)\|_{L^1(\Omega)} \le \|\omega_0\|_{L^1(\Omega)}$. Thus, if 
we define 
\[
  N(t) \,=\, \frac{\|\omega(t)\|_{L^2(\Omega)}^2}{\|\omega_0\|_{
  L^1(\Omega)}^2}~, \qquad \hbox{and}\quad \psi(x) \,=\, \min(x^2,x^3)~,
\]
we obtain the differential inequality $N'(t) \le -c \psi(N(t))$, 
for some constant $c > 0$. It follows that $\Psi(N(t)) \ge ct$ for all 
$t > 0$, where $\Psi : (0,\infty) \to (0,\infty)$ is the one-to-one
function defined by 
\[
  \Psi(x) \,=\, \int_x^\infty \frac{1}{\psi(y)}\dd y \,=\, 
  \begin{cases} ~\,\frac{1}{x} & x \ge 1~, \\[1mm] 
  \frac{x^2+1}{2x^2} & x < 1~,\end{cases} \qquad
  \Psi^{-1}(t) \,=\, \begin{cases} \quad\frac{1}{t} & t \le 1~, \\[1mm] 
  \frac{1}{\sqrt{2t-1}} & t > 1~.\end{cases}
\]
Since $\Psi$ is decreasing, we conclude that $N(t) \le \Psi^{-1}(ct)$ 
for all $t > 0$, hence
\[
   \|\omega(t)\|_{L^2(\Omega)}^2 \,\le\, \|\omega_0\|_{L^1(\Omega)}^2 
   \Psi^{-1}(ct) \,\le\, \|\omega_0\|_{L^1(\Omega)}^2 V(ct)^{-1}~,
  \qquad t > 0~,
\]
where $V(t) = \min(t,\sqrt{t})$. This shows that \eqref{LpLq} holds 
for $p = 1$, $q = 2$. 

To complete the proof, we use a classical duality argument. Fix $T > 0$ 
and let $w(x,t)$ be the solution of the adjoint equation
\begin{equation}\label{weq}
  \partial_t w - u\cdot \nabla w \,=\, \Delta w~, 
  \qquad x \in \Omega~, \quad 0 \le t \le T~,
\end{equation}
with initial data $w_0 \in L^1(\Omega)$. By construction the 
quantity $\int_\Omega \omega(x,t)w(x,T-t)\dd x$ is independent 
of time, so that
\[
  \int_\Omega \omega(x,T) w_0(x)\dd x \,=\, \int_\Omega \omega_0(x)
  w(x,T)\dd x~.
\]
Applying \eqref{LpLq} with $p = 1$, $q = 2$ to the adjoint equation 
\eqref{weq}, we thus obtain
\[
  \Bigl|\int_\Omega \omega(x,T) w_0(x)\dd x\Bigr| \,\le\, 
  \|\omega_0\|_{L^2(\Omega)} \|w(T)\|_{L^2(\Omega)} \,\le\, 
  C \|\omega_0\|_{L^2(\Omega)} \|w_0\|_{L^1(\Omega)} V(T)^{-1/2}~, 
\]
and since $w_0 \in L^1(\Omega)$ was arbitrary we conclude that 
$\|\omega(T)\|_{L^\infty(\Omega)} \le C \|\omega_0\|_{L^2(\Omega)} 
V(T)^{-1/2}$. This proves \eqref{LpLq} for $p = 2$, $q = \infty$.  

Finally, combining both estimates above, we obtain for any $t > 0$\:
\[
  \|\omega(t)\|_{L^\infty(\Omega)} \,\le\, \frac{C}{V(t/2)^{1/2}}
  \|\omega(t/2)\|_{L^2(\Omega)} \,\le\, \frac{C^2}{V(t/2)}
  \|\omega_0\|_{L^1(\Omega)}~,
\]
which proves \eqref{LpLq} for $p = 1$, $q = \infty$. \QED

\subsection{Bounds on the fundamental solution}\label{subfond}

The solution of \eqref{Veq3} with initial data $\omega_0$ can be 
represented as
\begin{equation}\label{omrep}
  \omega(x,t) \,=\, \int_\Omega \Gamma_u(x,y;t,0) \omega_0(y)\dd y~,
  \qquad x \in \Omega~, \quad t > 0~, 
\end{equation}
where $\Gamma_u(x,y;t,t_0)$ is the {\em fundamental solution} of the 
advection-diffusion equation \eqref{Veq3}. The strong maximum 
principle implies that $\Gamma_u(x,y;t,t_0) > 0$ whenever $t > t_0$, 
and it is also known that
\[
  \int_\Omega \Gamma_u(x,y;t,t_0)\dd y \,=\, 1~,\qquad \hbox{and}\quad
  \int_\Omega \Gamma_u(x,y;t,t_0)\dd x \,=\, 1~,
\]
for all $x,y \in \Omega$ and all $t > t_0$ (the last relation 
uses the assumption that $\div u = 0$). Finally, the semigroup 
property $\Gamma_u(x,y;t,t_0) = \int_\Omega \Gamma_u(x,z;t,s)
\Gamma_u(z,y;s,t_0)\dd z$ holds for all $x,y \in \Omega$ whenever 
$t > s > t_0$. We are interested in pointwise upper bounds on 
the fundamental solution $\Gamma_u$, in the spirit of Aronson 
\cite{Ar}. 

We assume henceforth that 
the velocity field $u(x,t)$ is uniformly bounded, and that the 
first component $u_1$ has zero vertical average\:
\begin{equation}\label{zeromean}
  \int_\T u_1(x_1,x_2,t)\dd x_2 \,=\, 0~, \qquad x_1 \in \R~, 
  \quad t \ge 0~.
\end{equation}
Under these hypotheses, we prove the following Gaussian upper bound 
on the fundamental solution. 

\begin{proposition}\label{FundSol}
Assume that $u$ is a divergence-free vector field satisfying 
\eqref{zeromean}, and such that $\sup_{t \ge 0}\|u_1(t)\|_{L^\infty(\Omega)} 
\le M$ for some $M \ge 0$. Then, for any $\lambda \in (0,1)$, there 
exists a constant $K_2 > 0$ (independent of $u$) such that  
\begin{equation}\label{Gaussbd}
  \Gamma_u(x,y;t,0) \,\le\, \frac{K_2}{V(t)}\, \exp\Bigl(-\lambda
  \frac{|x_1-y_1|^2}{4(1{+}M^2)t}\Bigr)~, \qquad x,y \in \Omega~, 
  \quad t > 0~,
\end{equation}
where $V(t) = \min(t,\sqrt{t})$. 
\end{proposition}

\begin{remark}\label{Fundstronger}
It is clear that \eqref{Gaussbd} implies estimate \eqref{LpLq}
for $p = 1$, $q = \infty$, and (as was already observed) the 
general case easily follows by interpolation. However, the proof of 
Proposition~\ref{FundSol} is substantially more complicated 
than that of Proposition~\ref{LpLqProp}.
\end{remark}

\noindent{\bf Proof.} We follow the approach of Fabes and Stroock
\cite{FS}. Let $\omega_0 : \Omega \to \R$ be continuous and 
compactly supported, and assume moreover that $\omega_0 \ge 0$ and 
$\omega_0 \not\equiv 0$. By the maximum principle, the solution of 
\eqref{Veq3} with initial data $\omega_0$ satisfies $\omega(x,t) > 0$ 
for all $x \in \Omega$ and all $t > 0$. Given any $\alpha \in \R$, we 
define
\begin{equation}\label{wdef}
  w(x,t) \,=\, e^{-\alpha x_1}\,\omega(x,t)~, \qquad x = (x_1,x_2) \in 
  \Omega~, \quad t \ge 0~.
\end{equation}
The new function $w(x,t)$ satisfies the modified equation
\begin{equation}\label{Weq}
  \partial_t w + u\cdot\nabla w + \alpha u_1 w \,=\, 
  \Delta w + 2\alpha \partial_1 w + \alpha^2 w~,
\end{equation}
where $u_1$ is the horizontal component of the velocity field $u$.
Since $u$ is divergence-free, a direct calculation shows that, 
for any positive integer $p \in \N^*$, 
\[
 \frac{\D}{\D t}\int_\Omega w(x,t)^{2p}\dd x \,=\, -2p(2p{-}1)\int_\Omega
 w^{2p-2}|\nabla w|^2 \dd x + 2p\alpha^2 \int_\Omega w^{2p} \dd x 
 - 2p\alpha \int_\Omega u_1w^{2p}\dd x~.
\]
As $p \ge 1$, we have
\[
  2p(2p{-}1)\int_\Omega w^{2p-2}|\nabla w|^2 \dd x \,\ge\, 2p^2 
  \int_\Omega w^{2p-2}|\nabla w|^2 \dd x \,=\, 2\int_\Omega 
  |\nabla w^p|^2 \dd x~.
\]
Moreover, it follows from that \eqref{zeromean} that $u_1 = \partial_2 v_1$ 
for some $v_1 : \Omega \to \R$ which satisfies the uniform bound
$\|v_1\|_{L^\infty(\Omega)} \le M/2$. Thus, integrating by parts, we obtain
\begin{align*}
  -2p\alpha \int_\Omega u_1 w^{2p}\dd x \,&=\, 
  -2p\alpha \int_\Omega \partial_2 v_1 (w^p)^2\dd x \,=\, 
  4p\alpha \int_\Omega v_1w^p\partial_2 w^p\dd x \\ 
  \,&\le\, 2p|\alpha| M \int_\Omega |w^p||\partial_2 w^p|\dd x \,\le\, 
  \int_\Omega \Bigl(|\nabla w^p|^2 + p^2 \alpha^2 M^2w^{2p}\Bigr)\dd x~.
\end{align*}
Combining these estimates, we arrive at
\begin{equation}\label{Fund1}
 \frac{\D}{\D t}\int_\Omega w(x,t)^{2p}\dd x \,\le\,  -\int_\Omega
 |\nabla w(x,t)^p|^2 \dd x + \alpha^2(2p + p^2M^2)\int_\Omega 
 w(x,t)^{2p} \dd x~.
\end{equation}

To simplify the notations we define, for all $p \ge 1$,  
\[
  w_p(t) \,=\, \|w(\cdot,t)\|_{L^p(\Omega)} \,=\, \Bigl(\int_\Omega 
  |w(x,t)|^p\dd x\Bigr)^{1/p}~, \qquad t \ge 0~.
\]
Applying Nash's inequality \eqref{psiNashbd} to the function 
$f = w(\cdot,t)^p > 0$, we obtain the lower bound
\[
  \int_\Omega|\nabla w(x,t)^p|^2\dd x \,\ge\, C w_{2p}(t)^{2p}
  \min\left\{\frac{w_{2p}(t)^{2p}}{w_p(t)^{2p}}\,,\,
  \frac{w_{2p}(t)^{4p}}{w_p(t)^{4p}}\right\}~,
\]
for some universal constant $C > 0$. Thus it follows from 
\eqref{Fund1} that
\begin{equation}\label{Fund2}
  w_{2p}'(t) \,\le\, -\frac{C}{2p} \min_{\beta = 2,4}\biggl\{\Bigl(
  \frac{w_{2p}(t)}{w_p(t)} \Bigr)^{\beta p}\biggr\}w_{2p}(t) + 
  \alpha^2 \Bigl(1 + \frac{p}{2}M^2\Bigr)w_{2p}(t)~, \qquad t > 0~.
\end{equation}

\begin{remark}\label{FSrem}
In \cite[Section~1]{FS} the authors obtain a differential 
inequality of the form \eqref{Fund2} with $\beta = 4/n$, 
where $n \in \N^*$ is the space dimension. Here we have 
a combination of $\beta = 4$ ($n = 1$) and $\beta = 2$ 
($n = 2$) because we consider the cylinder $\Omega = \R \times \T$.
This makes the inequalities \eqref{Fund2} less homogeneous
and more cumbersome to integrate. 
\end{remark}

\begin{lemma}\label{intinlem}
Assume that inequalities \eqref{Fund2} hold for all $p \in \{1\} 
\cup S$, where $S = \{2^k \,|\, k \in \N^*\}$. Then for any 
$\epsilon > 0$ there exists a constant $C_\epsilon > 0$ such that
\begin{equation}\label{intineq}
  w_p(t) \,\le\, C_\epsilon\,\frac{e^{(1+\epsilon)\alpha^2(1+M^2)t}}{
  V(t)^{\frac{p-2}{2p}}}\,w_2(0)~, \qquad t > 0~, \qquad p \in S~,
\end{equation}
where $V(t) = \min(t,\sqrt{t})$. 
\end{lemma}

The proof of Lemma~\ref{intinlem} being somewhat technical, we 
postpone it to Section~\ref{secApp}, and assuming that \eqref{intineq}
holds we now conclude the proof of Proposition~\ref{FundSol}. 
Taking the limit $p \to \infty$ in \eqref{intineq}, we obtain
\[
  \|w(t)\|_{L^\infty(\Omega)} \,\le\, \frac{C_\epsilon}{V(t)^{1/2}}
  \,\,e^{(1+\epsilon)\alpha^2(1+M^2)t}\,\|w(0)\|_{L^2(\Omega)}~, 
  \qquad t > 0~.
\]
As in the proof of Proposition~\ref{LpLqProp}, a duality argument
gives the same bound for $\|w(t)\|_{L^2(\Omega)}$ in terms of 
$\|w_0\|_{L^1(\Omega)}$, so altogether we obtain
\begin{equation}\label{Fund3}
  \|w(t)\|_{L^\infty(\Omega)} \,\le\, \frac{\tilde C_\epsilon}{V(t)}
  \,\,e^{(1+\epsilon)\alpha^2(1+M^2)t}\,\|w(0)\|_{L^1(\Omega)}~, 
  \qquad t > 0~.
\end{equation}
Finally, we return to the original equation \eqref{Veq3}. 
If we take a sequence of initial data $\omega_0$ approaching
a Dirac mass at some point $y \in \Omega$, the corresponding
solutions $\omega(x,t)$ converge by definition to the fundamental 
solution $\Gamma_u(x,y;t,0)$. In view of \eqref{wdef}, estimate 
\eqref{Fund3} then implies
\begin{equation}\label{Fund4}
  \Gamma_u(x,y;t,0) \,\le\, \frac{\tilde C_\epsilon}{V(t)}
  \,\,e^{(1+\epsilon)\alpha^2(1+M^2)t}\,e^{\alpha(x_1-y_1)}~,
\end{equation}
for all $x,y \in \Omega$, all $t > 0$, and all $\alpha \in \R$. 
If we now choose
\[
  \alpha \,=\, -\frac{x_1-y_1}{2(1{+}\epsilon)(1{+}M^2)t}~, 
\]
we obtain from \eqref{Fund4}
\begin{equation}\label{Fund5}
  \Gamma_u(x,y;t,0) \,\le\, \frac{\tilde C_\epsilon}{V(t)}\,
  \exp\Bigl(-\frac{|x_1-y_1|^2}{4(1{+}\epsilon)(1{+}M^2)t}\Bigr)~, 
  \qquad x,y \in \Omega~, \quad t > 0~.
\end{equation}
This gives \eqref{Gaussbd} if $\epsilon > 0$ is taken sufficiently 
small. \QED

\begin{remark}\label{Gaussrem2}
It does not seem possible to obtain estimate \eqref{Fund5}
for all times using the simpler argument given in the 
proof of Proposition~\ref{LpLqProp}. The reason is that, 
when $\alpha \neq 0$, we do not have a good a priori estimate
on the $L^1$ norm of $w(x,t)$. The best we can deduce from 
\eqref{Weq} is
\[
  \|w(t)\|_{L^1(\Omega)} \,\le\, \|w_0\|_{L^1(\Omega)}\,
  e^{(\alpha^2 + |\alpha|M)t}~, \qquad t > 0~,
\]
which does not take into account the crucial assumption 
\eqref{zeromean}, and therefore cannot be used to derive 
estimate \eqref{Gaussbd} for large times. 
\end{remark}

\section{Weighted energy estimates}\label{secEner}

We now begin the actual proof of Theorem~\ref{main}. In what 
follows, we always assume that $u(x,t)$ is a solution of the 
Navier-Stokes equations \eqref{NSeq2} in $\Omega = \R 
\times \T$ satisfying \eqref{zeromean}, with associated vorticity 
$\omega = \partial_1 u_2 - \partial_2 u_1$ and pressure $p$ given 
by \eqref{pdef}. As in Theorem~\ref{Cauchy}, we suppose that the 
initial velocity $u_0$ is divergence-free, and that $\omega_0 
= \curl u_0 \in L^\infty(\Omega)$. We denote $M = 
\|\omega_0\|_{L^\infty}$. It then follows from the maximum
principle that the vorticity $\omega(x,t)$, which solves
\eqref{Veq2}, satisfies $\|\omega(t)\|_{L^\infty} \le M$ for all 
$t \ge 0$. 

\subsection{Energy density, energy flux, energy dissipation}
\label{subeee}

As in \cite{GS3}, our approach relies on a careful study of 
the local energy dissipation in the system. For any $x_1 \in \R$ 
and $t \ge 0$, we define
\begin{align}\label{edef}
  e(x_1,t) \,&=\, \frac12\int_\T |u(x_1,x_2,t)|^2 \dd x_2
    + \frac{M^2}{2}~, \\ \label{hdef}
  h(x_1,t) \,&=\, \int_\T \Bigl(p(x_1,x_2,t) + \frac12|u(x_1,x_2,t)|^2
    \Bigr)u_1(x_1,x_2,t) \dd x_2~, \\ \label{ddef}
  d(x_1,t) \,&=\, \int_\T |\nabla u(x_1,x_2,t)|^2 \dd x_2~,
\end{align}
as well as $f(x_1,t) = \partial_1e(x_1,t) - h(x_1,t)$. The quantities
$e,f,d$ will be referred to as the energy density, the energy flux, 
and the energy dissipation rate, respectively. It is clear that 
$e(x_1,t) \ge 0$ and $d(x_1,t) \ge 0$. Moreover, a direct calculation
shows that the following local energy dissipation law holds\: 
\begin{equation}\label{led}
  \partial_t e(x_1,t) \,=\, \partial_1 f(x_1,t) - d(x_1,t)~, 
  \qquad x_1 \in \R~, \qquad t > 0~.
\end{equation}
Finally, the initial energy density is uniformly bounded, and 
we have
\begin{equation}\label{e0bound}
  e_*(0) \,=\, \sup_{x_1 \in \R}e(x_1,0) \,\le\, \frac12 
  \|u_0\|_{L^\infty}^2 + \frac{M^2}{2}~.
\end{equation}

The reason for including the constant $M^2/2$ in the definition
\eqref{edef} of the energy density will become clear in the proof of
the following lemma, which provides a crucial estimate of the energy
flux in terms of the energy density and the energy dissipation.

\begin{lemma}\label{fluxbound}
There exists a constant $C_3 > 0$ such that 
\begin{equation}\label{fluxbd}
  |f|^2 \,\le\, C_3(1+M)^2e d~, \qquad |\partial_1 e|^2 \,\le\, 
  2e d~.
\end{equation}
\end{lemma}

\noindent{\bf Proof.} We first estimate the inviscid flux $h$ 
defined by \eqref{hdef}. Since $\langle u_1 \rangle = 0$, the 
Poincar\'e-Wirtinger inequality implies that
\[
  \int_\T |u_1(x_1,x_2,t)|^2 \dd x_2 \,\le\, \frac{1}{4\pi^2} 
  \int_\T |\partial_2 u_1(x_1,x_2,t)|^2 \dd x_2 \,\le\, 
  \frac{d(x_1,t)}{4\pi^2}~. 
\]
Using \eqref{pres} and H\"older's inequality, we thus find 
\begin{equation}\label{flux1}
  \Bigl|\int_\T p u_1 \dd x_2\Bigr| \,\le\, C_2 M^2 \int_\T |u_1| 
  \dd x_2 \,\le\, \frac{C_2M^2}{2\pi}\,d^{1/2} \,\le\, 
  \frac{C_2M}{2\pi}\,(2ed)^{1/2}~, 
\end{equation}
where in the last inequality we used the fact that $e \ge M^2/2$. 
On the other hand, we have $u_1 = \widehat u_1 = \partial_2 v_1$, where 
$\|v_1\|_{L^\infty} \le \frac12 \|\widehat u_1\|_{L^\infty} \le C_1 M$ 
by \eqref{BSres}. Therefore
\[
  \frac12 \int_\T |u|^2 u_1 \dd x_2 \,=\, -\int_\T (u\cdot 
  \partial_2 u)v_1 \dd x_2~,
\]
and using H\"older's inequality again we obtain
\begin{equation}\label{flux2}
  \left|\frac12 \int_\T |u|^2 u_1 \dd x_2\right| \,\le\, 
  \|v_1\|_{L^\infty} \!\int_\T |u||\partial_2 u|\dd x_2 
  \,\le\, C_1 M (2ed)^{1/2}~.
\end{equation}
Inequalities \eqref{flux1} and \eqref{flux2} imply that 
$|h|^2 \le CM^2 ed$ for some $C > 0$. On the other hand, we have 
$\partial_1 e = \int_\T u\cdot\partial_1 u\dd x_2$, hence 
$|\partial_1 e|^2 \le 2ed$ by H\"older's inequality.  Since
$f = \partial_1e - h$, this proves \eqref{fluxbd}. \QED

\subsection{Localized energy estimate}\label{subener}

In what follows we denote $\beta = C_3(1+M)^2$, where $C_3 > 0$ 
is as in \eqref{fluxbd}. By Lemma~\ref{fluxbound}, the energy flux
satisfies
\begin{equation}\label{fluxbd2}
  |f(x_1,t)|^2 \,\le\,  \beta e(x_1,t) d(x_1,t)~, \qquad 
  x_1 \in \R~, \quad t > 0~.
\end{equation}
Our goal is to control the solution \eqref{NSeq2} using localized 
energy estimates. Given $\rho > 0$, we introduce the localization 
function $\chi_\rho(x_1) = \exp(-\rho|x_1|)$, and we define
\begin{equation}\label{EDdef}
  E_\rho(t) \,=\, \int_\R \chi_\rho(x_1) e(x_1,t)\dd x_1~, \qquad 
  D_\rho(t) \,=\, \int_\R \chi_\rho(x_1) d(x_1,t)\dd x_1~,
\end{equation}
for all $t \ge 0$. We then have the following estimate on the 
localized energy $E_\rho(t)$\:

\begin{proposition}\label{Elem}
Fix $T > 0$, and let $\rho = 1/\sqrt{\beta T}$ where $\beta > 0$ 
is as in \eqref{fluxbd2}. Then
\begin{equation}\label{Ebound}
  E_\rho(T) + \frac12 \int_0^T D_\rho(t)\dd t \,\le\, 
  4 e_*(0)\sqrt{\beta T}~,
\end{equation}
where $e_*(0)$ is given by \eqref{e0bound}. 
\end{proposition}

\noindent{\bf Proof.} Differentiating $E_\rho(t)$ with respect to 
time and using \eqref{led}, we obtain
\begin{equation}\label{Ediff1}
  E_\rho'(t) \,=\, \int_\R \chi_\rho \partial_t e\dd x_1 \,=\, 
  \int_\R \chi_\rho (\partial_1 f - d)\dd x_1 \,=\, 
  - \int_\R \Bigl(\chi_\rho' f + \chi_\rho d\Bigr)\dd x_1~.
\end{equation}
Since $|\chi_\rho'| \le \rho \chi_\rho$, it follows from
\eqref{fluxbd2} that
\[
  \Bigl|\int_\R \chi_\rho' f \dd x_1\Bigr| \,\le\, 
  \rho \int_\R \chi_\rho (\beta e d)^{1/2}\dd x_1 \,\le\, 
  \frac12 \int_\R \chi_\rho d\dd x_1 + \frac{\rho^2\beta}{2} 
  \int_\R \chi_\rho e\dd x_1~. 
\]
Thus \eqref{Ediff1} leads to the differential inequality 
$E_\rho'(t) +  \frac12 D_\rho(t) \le \frac12\,\rho^2 \beta 
E_\rho(t)$, which can be integrated using Gronwall's lemma
to give
\[
  E_\rho(T) + \frac12 \int_0^T D_\rho(t)\dd t \,\le\, 
  E_\rho(0)\exp\Bigl(\frac12\,\rho^2 \beta T\Bigr)~.
\]
Since $E_\rho(0) \le e_*(0)\int_\R \chi_\rho(x_1)\dd x_1 = 
2e_*(0)/\rho$, choosing $\rho = (\beta T)^{-1/2}$ yields
the desired result. \QED

\begin{remark}\label{compare}
Together with \eqref{led}, Lemma~ \ref{fluxbound} implies
that the Navier-Stokes equations in the domain $\Omega$ 
define a one-dimensional "extended dissipative system", in 
the sense of \cite{GS2}. This point of view was thoroughly 
exploited in \cite{GS3}, where results similar to 
Proposition~\ref{Elem} were obtained using a slightly 
different approach. In particular, one can verify that 
estimate \eqref{uintbound} with $R = \sqrt{\nu T}$ is 
equivalent to \eqref{Ebound}. 
\end{remark}

\subsection{Localized enstrophy estimate}\label{subens}

We now perform a similar analysis at the level of the vorticity 
equation \eqref{Veq2}. In analogy with \eqref{edef}--\eqref{ddef} 
we define, for all $x_1 \in \R$ and all $t \ge 0$,
\begin{align}\label{eedef}
  \varepsilon(x_1,t) \,&=\, \frac12\int_\T |\omega(x_1,x_2,t)|^2 \dd x_2
   ~, \\ \label{hhdef}
  \zeta(x_1,t) \,&=\, \frac12 \int_\T \omega(x_1,x_2,t)^2
    u_1(x_1,x_2,t) \dd x_2~, \\ \label{dddef}
  \delta(x_1,t) \,&=\, \int_\T |\nabla \omega(x_1,x_2,t)|^2 \dd x_2~,
\end{align}
as well as $\phi(x_1,t) = \partial_1\varepsilon(x_1,t) - \zeta(x_1,t)$. 
Using \eqref{Veq2} we obtain the local enstrophy 
dissipation law
\begin{equation}\label{tildeled}
  \partial_t \varepsilon(x_1,t) \,=\, \partial_1 \phi(x_1,t) - 
  \delta(x_1,t)~, \qquad x_1 \in \R~, \qquad t > 0~,
\end{equation}
and we have the analog of Lemma~\ref{fluxbound}\:

\begin{lemma}\label{tfluxbound}
There exists a constant $C_4 > 0$ such that 
\begin{equation}\label{tfluxbd}
  |\phi|^2 \,\le\, C_4(1+M)^2\varepsilon \delta~, \qquad 
  |\partial_1 \varepsilon|^2 \,\le\, 2\varepsilon \delta~.
\end{equation}
\end{lemma}

\noindent{\bf Proof.} We proceed as in the proof of 
Lemma~\ref{fluxbound}. As $u_1 = \partial_2 v_1$ with 
$\|v_1\|_{L^\infty} \le C_1 M$, we have
\[
  |\zeta| \,=\, \left|\int_\T \omega(\partial_2\omega) v_1 \dd x_2
  \right| \,\le\, \|v_1\|_{L^\infty} \!\int_\T |\omega||\partial_2 
  \omega|\dd x_2 \,\le\, C_1 M (2\varepsilon\delta)^{1/2}~.
\]
Since $\phi = \partial_1\varepsilon - \zeta$ and $|\partial_1 
\varepsilon|^2 \le 2\varepsilon\delta$ by H\"older's inequality, 
we obtain \eqref{tfluxbd}. \QED

\bigskip
As in \eqref{EDdef} we define the localized enstrophy and
the corresponding dissipation by 
\begin{equation}\label{tEDdef}
  \EE_\rho(t) \,=\, \int_\R \chi_\rho(x_1) \varepsilon(x_1,t)
  \dd x_1~, \qquad \DD_\rho(t) \,=\, \int_\R \chi_\rho(x_1) 
  \delta(x_1,t)\dd x_1~.
\end{equation}
We then have the following estimate\:

\begin{proposition}\label{tElem}
There exists $C_5 > 0$ such that, if $T > 0$ and $\rho = 
1/\sqrt{\beta T}$, then 
\begin{equation}\label{tEbound}
  \EE_\rho(T) + \frac12 \int_{T/2}^T \DD_\rho(t)\dd t 
  \,\le\, C_5(1+M)e_*(0)\frac{1}{\sqrt{T}}~.
\end{equation} 
\end{proposition}

\noindent{\bf Proof.} Proceeding as in the proof of 
Proposition~\ref{Elem}, we obtain for $\EE_\rho(t)$ 
the differential inequality
\begin{equation}\label{tE1}
  \EE_\rho'(t) + \frac12 \DD_\rho(t) \,\le\, 
  \frac{C_4}{2}(1+M)^2\rho^2 \EE_\rho(t)~, \qquad 
  t > 0~.
\end{equation}
Since $\omega^2 \le 2|\nabla u|^2$, we have $\varepsilon(x_1,t) 
\le d(x_1,t)$, hence $\EE_\rho(t) \le D_\rho(t)$ for all $t \ge 0$. 
Therefore \eqref{Ebound} implies
\begin{equation}\label{tE2}
  \int_0^T \EE_\rho(t)\dd t \,\le\, \int_0^T D_\rho(t)\dd t \,\le\, 
  8 e_*(0)\sqrt{\beta T}~.
\end{equation}
In particular, there exists $t_0 \in [0,T/2]$ such that $\EE_\rho(t_0) 
\le 16 e_*(0)\sqrt{\beta/T}$. Integrating \eqref{tE1}
between $t_0$ and $T$ and using \eqref{tE2}, we thus obtain
\begin{equation}\label{tE3}
  \EE_\rho(T) + \frac12 \int_{t_0}^T \DD_\rho(t) \dd t
  \,\le\, \EE_\rho(t_0) + 4C_4(1+M)^2\rho^2 e_*(0)
  \sqrt{\beta T}~.
\end{equation}
This gives the desired result since $t_0 \le T/2$, $\rho = 
1/\sqrt{\beta T}$, and $\beta = C_3(1+M)^2$. \QED

\section{Uniform estimates for the velocity and the 
vorticity}\label{secUnif}

Combining the weighted energy estimates of the previous section with
the bounds on the fundamental solution of the vorticity equation
obtained in Section~\ref{secVort}, we are now able to prove assertions
1) and 2) in Theorem~\ref{main}.  We keep the same notations as in
Section~\ref{secEner}. In particular, $u(x,t)$ is a solution of
\eqref{NSeq2} with bounded initial velocity $u_0$ and vorticity
$\omega_0$, which satisfies \eqref{zeromean}, and we denote $M =
\|\omega_0\|_{L^\infty}$.

\subsection{Uniform decay of the vorticity}

\begin{proposition}\label{Unifom}
There exists a constant $C_6 > 0$ such that
\begin{equation}\label{omdecay}
  \|\omega(t)\|_{L^\infty}^2 \,\le\, C_6(1+M)e_*(0)\frac{1}{\sqrt{t}}~,
  \qquad t > 0~,
\end{equation}
where $M = \|\omega_0\|_{L^\infty}$ and $e_*(0)$ is defined in 
\eqref{e0bound}. 
\end{proposition}

\noindent{\bf Proof.} Since $\|\omega(t)\|_{L^\infty} \le M$ for all
$t \ge 0$ and $e_*(0) \ge M^2/2$, it is clear that \eqref{omdecay}
holds with $C_6 = 2\sqrt{2}$ whenever $t \le 2(1+M)^2$. Thus we assume
henceforth that $t \ge 2(1+M)^2$, and given such a time $t$ we denote $T
= t-1 \ge t/2 \ge 1$. We also define $A = \sqrt{\beta T}$, where 
$\beta = C_3(1+M)^2$ is as in \eqref{fluxbd2}. Using the fundamental 
solution $\Gamma_u$ introduced in Section~\ref{secVort}, we decompose
\[
  \omega(x,t) \,=\, \int_{\Omega_1} \Gamma_u(x,y;t,T) \omega(y,T)\dd y
 + \int_{\Omega_2} \Gamma_u(x,y;t,T) \omega(y,T)\dd y \,=\, 
 \omega_1(x,t) + \omega_2(x,t)~,
\]
where $\Omega_1 = \{x \in \Omega\,|\, |x_1| \le A\}$ and $\Omega_2 =
\{x \in \Omega\,|\, |x_1| > A\}$. In view of Propositions~\ref{LpLqProp} 
and \ref{tElem}, we have
\begin{equation}\label{om1bd}
  \sup_{x \in \Omega}|\omega_1(x,t)|^2 \,\le\, K_1^2 \int_{\Omega_1}
  |\omega(y,T)|^2 \dd y \,\le\, K_1^2 (2eC_5) (1+M)e_*(0)
  \frac{1}{\sqrt{T}}~. 
\end{equation}
To bound $\omega_2$, we use Proposition~\ref{FundSol} which gives, 
for any $\lambda \in (0,1)$, 
\begin{equation}\label{om2bd1}
  |\omega_2(x,t)| \,\le\, K_2 \int_{\Omega_2} \exp\Bigl(-\lambda
  \frac{|x_1-y_1|^2}{4(1{+}M)^2}\Bigr)|\omega(y,T)|\dd y~,
  \qquad x \in \Omega~.
\end{equation}
In particular, if $|x_1| \le A/2$, we have $|x_1 - y_1| \ge A/2$ whenever 
$y \in \Omega_2$, hence using the a priori estimate $|\omega(y,T)| \le M$
we find
\[
  |\omega_2(x,t)| \,\le\, 2 K_2M \int_{A/2}^\infty \exp\Bigl(-
  \frac{\lambda z^2}{4(1{+}M)^2}\Bigr)\dd z \,\le\, 8 K_2 M \frac{(1{+}M)^2}{
  A\lambda}\,\exp\Bigl(-\frac{\lambda A^2}{16(1{+}M)^2}\Bigr)~.
\]
Since $A^2 = \beta T$ with $\beta = C_3(1+M)^2$ and $T \ge (1{+}M)^2$, 
we conclude that 
\begin{equation}\label{om2bd2}
  \sup_{|x_1| \le A/2}|\omega_2(x,t)| \,\le\, CM(1{+}M)\frac{1}{\sqrt{T}}
  \,e^{-\lambda C_3 T/16} \,\le\, \frac{CM}{\sqrt{T}}~, 
\end{equation}
for some $C > 0$. Combining \eqref{om1bd}, \eqref{om2bd2} and using
the fact that $T \ge t/2 \ge 1$, we obtain
\begin{equation}\label{om12bd}
  \sup_{|x_1| \le A/2}|\omega(x,t)|^2 \,\le\, C_6(1+M)e_*(0)
  \frac{1}{\sqrt{t}}~, \qquad t \ge 2(1+M)^2~,
\end{equation}
for some $C_6 > 0$. Now, it is clear that estimate \eqref{om12bd} still 
holds if we replace the vorticity $\omega(x_1,x_2,t)$ by $\omega(x_1{-}a,
x_2,t)$ for any $a \in \R$, because equations \eqref{NSeq2}, \eqref{Veq2}
are translation invariant in the horizontal direction, and the right-hand side 
of \eqref{om12bd} involves only translation invariant quantities. Thus 
in \eqref{om12bd} we can take the supremum over all $x \in \Omega$, 
and this proves that \eqref{omdecay} holds for all $t \ge 2(1+M)^2$. \QED

\subsection{Uniform bound on the velocity field}

\begin{proposition}\label{Unifu}
There exists a constant $C_7 > 0$ such that
\begin{equation}\label{uunif}
  \|u(t)\|_{L^\infty} \,\le\, C_7\Bigl(\|u_0\|_{L^\infty} + M + 
  (1+M)e_*(0)\Bigr)~, \qquad t \ge 0~,
\end{equation}
where $M = \|\omega_0\|_{L^\infty}$ and $e_*(0)$ is defined in 
\eqref{e0bound}. 
\end{proposition}

\noindent{\bf Proof.} If $u(x,t)$ is decomposed as in \eqref{udecomp}, 
we already know that $\|\widehat u(t)\|_{L^\infty} \le 2C_1 \|\omega(t)
\|_{L^\infty} \,\le\, 2C_1 M$ for all $t \ge 0$. Thus it remains to 
bound the mean vertical flux $m(x_1,t)$. The integral equation 
corresponding to \eqref{meq} is
\[
  m(t) \,=\, S_1(t)m(0) - \int_0^t \partial_1 S_1(t-s)
  \langle \widehat u_1(s)\widehat u_2(s)\rangle\dd s~,
\]
where $S_1(t) = e^{t \partial_1^2}$ is heat semigroup on $\R$. 
By Proposition~\ref{Unifom}, we have
\[
  \|\widehat u(t)\|_{L^\infty}^2 \,\le\, 4C_1^2 \|\omega(t)\|_{L^\infty}^2 
  \,\le\, 4C_1^2 C_6(1+M)e_*(0)\frac{1}{\sqrt{t}}~, \qquad t > 0~,
\]
hence
\[
  \|m(t)\|_{L^\infty} \,\le\, \|m(0)\|_{L^\infty} + \int_0^t 
  \frac{\|\widehat u(s)\|_{L^\infty}^2}{\sqrt{\pi(t-s)}}
  \dd s \,\le\, \|m(0)\|_{L^\infty} + 4\sqrt{\pi}C_1^2 C_6(1+M)e_*(0)~.
\]
We conclude that, for all $t \ge 0$, 
\[
  \|u(t)\|_{L^\infty} \,\le\, \|\widehat u(t)\|_{L^\infty} + \|m(t)\|_{L^\infty} 
  \,\le\, 2C_1 M + \|u_0\|_{L^\infty} + 4\sqrt{\pi}C_1^2 C_6(1+M)e_*(0)~, 
\]
and \eqref{uunif} follows. \QED 

\begin{remark}\label{back}
According to \eqref{adim}, to translate our results back into the 
original variables we have to replace  $\|u(t)\|_{L^\infty}$ by 
$L\|u(t)\|_{L^\infty}/\nu$, $\|\omega(t)\|_{L^\infty}$ by 
$L^2\|\omega(t)\|_{L^\infty}/\nu$, and $t$ by $\nu t/L^2$. 
Thus $\|u_0\|_{L^\infty}$ is replaced by $R_u$, and $M = 
\|\omega_0\|_{L^\infty}$ by $R_\omega$. Since $e_*(0) \le
\frac12 M^2 + \frac12 \|u_0\|_{L^\infty}^2$, we see that 
\eqref{main1}, \eqref{main2} follow from \eqref{uunif}, 
\eqref{omdecay} respectively. 
\end{remark}

\section{Exponential decay in the laminar regime}\label{secExp}

Finally we prove assertion 3) in Theorem~\ref{main}. As is clear 
from Proposition~\ref{Unifom}, any solution of \eqref{NSeq2}, 
\eqref{pdef} with bounded initial data satisfies $\|\omega(t)\|_{L^\infty} 
< 4\pi^2$ for $t$ sufficiently large. In this section, we assume that 
the initial vorticity $\omega_0$ is small enough so that
\begin{equation}\label{kappa}
  \kappa \,:=\, \frac{\|\omega_0\|_{L^\infty}}{4\pi^2} \,<\, 1~.
\end{equation}
If the solution $u(x,t)$ is decomposed as in \eqref{udecomp}, 
we define, in analogy with \eqref{edef}--\eqref{ddef}, 
\begin{align} \label{hatedef}
  \widehat e(x_1,t) \,&=\, \frac{1}{2}\int_\T |\widehat 
  u(x_1,x_2,t)|^2 \dd x_2~, \\ \label{hathdef}
  \widehat h(x_1,t) \,&=\, \int_\T \Bigl(p + \frac12
  |\widehat u(x_1,x_2,t)|^2\Bigr)\widehat u_1(x_1,x_2,t) 
  \dd x_2~, \\ \label{hatddef} \widehat d(x_1,t) \,&=\, 
  \int_\T |\nabla \widehat u(x_1,x_2,t)|^2 \dd x_2~,
\end{align}
as well as $\widehat f(x_1,t) = \partial_1\widehat e(x_1,t) - 
\widehat h(x_1,t)$. Using \eqref{hatu1eq}, \eqref{hatu2eq}, 
it is not difficult to establish the modified energy dissipation law
\begin{equation}\label{hateds}
  \partial_t \widehat e(x_1,t) \,=\, \partial_1 \widehat f(x_1,t) - 
  \widehat d(x_1,t) - \widehat g(x_1,t)~, \qquad x_1 \in \R~, 
  \qquad t > 0~,
\end{equation}
where $\widehat g \,=\, (\partial_1 m) \langle \widehat u_1 \widehat 
u_2\rangle$. As in Lemma~\ref{fluxbound}, we then have

\begin{lemma}\label{hatfluxbound}
There exists a constant $C_8 > 0$ such that 
\begin{equation}\label{hatfluxbd}
  \widehat e \,\le\, \frac{\widehat d}{8\pi^2}~, \qquad
  |\widehat f|^2 \,\le\, C_8\kappa^2\,\widehat d~, \qquad
  |\widehat g| \,\le\, \kappa \widehat d~.
\end{equation}
\end{lemma}

\noindent{\bf Proof.} The first and the last estimate in 
\eqref{hatfluxbd} follow immediately from the Poincar\'e-Wirtinger
inequality, if we observe in addition that  $|\partial_1 m| \le 
|\langle \omega\rangle| \le \|\omega_0\|_{L^\infty}$. To bound 
$\widehat h$, we proceed as in \eqref{flux1} and \eqref{flux2}. 
We find
\[
  \Bigl|\int_\T p \,\widehat u_1 \dd x_2\Bigr| \,\le\, 
  C_2 \|\omega_0\|_{L^\infty}^2 \int_\T |\widehat u_1| \dd x_2 \,\le\, 
  C \kappa^2\,\widehat d\,{}^{1/2}~, 
\]
and
\[
  \left|\frac12 \int_\T |\widehat u|^2 \widehat u_1 \dd x_2\right| 
  \,\le\, \|v_1\|_{L^\infty} \!\int_\T |\widehat u||\partial_2 \widehat 
  u|\dd x_2 \,\le\, C \kappa^2\,\widehat d\,{}^{1/2}~, 
\]
In the last inequality, we used the fact that both $\|v_1\|_{L^\infty}$
and $\|\widehat u\|_{L^\infty}$ are bounded by $C_1 \|\omega_0\|_{L^\infty}$. 
Thus we have $|\widehat h|^2 \le C\kappa^4 \widehat d$, and we 
also know that $|\partial_1 \widehat e|^2 \le  2 \widehat e\widehat d 
\le C\kappa^2 \widehat d$. Combining these estimates and using the 
assumption that $\kappa < 1$, we obtain $|\widehat f|^2 \,\le\, 
C_8\kappa^2\,\widehat d$ for some $C_8 > 0$. \QED

\begin{proposition}\label{hatudecay}
If the initial data satisfy \eqref{kappa}, then for any $\epsilon > 0$ 
there exists a constant $C_9 > 0$ such that
\begin{equation}\label{hatubd}
  \sup_{a \in \R}  \int_{a-1}^{a+1}\widehat e(x_1,t)\dd x_1 \,\le\, 
  \frac{C_9 \kappa^2}{1-\kappa}\,e^{-\gamma t/2}~, \qquad t \ge 0~,
\end{equation}
where $\gamma = 8\pi^2(1-\epsilon)(1-\kappa)$. 
\end{proposition}

\noindent{\bf Proof.} Following the approach developed in 
\cite{GS2,GS3}, we first establish a differential inequality 
for the energy density $\widehat e$ defined in \eqref{hatedef}. 
Using \eqref{hateds} and \eqref{hatfluxbd}, we find
\begin{equation}\label{hatproto}
  \partial_t \widehat e(x_1,t) \,\le\, \partial_1 \widehat f(x_1,t) 
  - (1-\kappa)\widehat d(x_1,t) \,\le\, \partial_1 \widehat f(x_1,t) 
  - \eta \widehat f(x_1,t)^2 - \gamma \widehat e(x_1,t)~,
\end{equation}
where $\eta,\gamma > 0$ satisfy $C_8\kappa^2 \eta + \gamma/(8\pi^2) 
= 1-\kappa$. For definiteness, we take $\epsilon \in (0,1)$ and choose
\begin{equation}\label{bgdef}
  \eta \,=\, \frac{\epsilon(1-\kappa)}{C_8\kappa^2}~, \qquad 
  \quad \gamma \,=\, 8\pi^2(1-\epsilon)(1-\kappa)~.
\end{equation}
Inequality \eqref{hatproto} can be written in the equivalent form
\begin{equation}\label{hatediff}
  \partial_t \Bigl(\widehat e(x_1,t)\,e^{\gamma t}\Bigr) \,\le\, e^{\gamma t}
  \Bigl(\partial_1 \widehat f(x_1,t) - \eta \widehat f(x_1,t)^2\Bigr)~, 
  \qquad x_1 \in \R~, \quad t> 0~.
\end{equation}
To exploit \eqref{hatediff}, we fix $T > 0$ and we introduce the 
integrated flux
\[
  F(x_1) \,=\, \int_0^T \widehat f(x_1,t)e^{\gamma t}\dd t~, \qquad 
  x_1 \in \R~,
\]
which can be estimated as follows\:
\[
  F(x_1)^2 \,\le\, \biggl(\int_0^T e^{\gamma t}\dd t\biggr) \biggl(\int_0^T 
  \widehat f(x_1,t)^2 e^{\gamma t}\dd t\biggr) \,\le\, 
  \frac{e^{\gamma T}}{\gamma}\biggl(\int_0^T \widehat f(x_1,t)^2 
  e^{\gamma t}\dd t\biggr)~.
\]
Integrating both sides of \eqref{hatediff} over $t \in [0,T]$, we 
thus obtain
\begin{equation}\label{Fdiff1}
  \widehat e(x_1,T)\,e^{\gamma T} - \widehat e(x_1,0) \,\le\, 
  F'(x_1) - \eta\gamma e^{-\gamma T}F(x_1)^2~, \qquad x_1 \in \R~. 
\end{equation}
In particular, we see that the integrated flux $F(x_1)$ satisfies 
the differential inequality
\begin{equation}\label{Fdiff2}
  F'(x_1) \,\ge\, - \widehat e_*(0) + \eta\gamma e^{-\gamma T}F(x_1)^2~,
  \qquad x_1 \in \R~,
\end{equation}
where $\widehat e_*(0) = \sup_{x_1 \in \R}\widehat e(x_1,0) \le 2
C_1^2 \|\omega_0\|_{L^\infty}^2$. Now, it is easy to verify that any
solution of \eqref{Fdiff2} that is globally defined on $\R$ 
necessarily satisfies
\begin{equation}\label{Fbound}
  F(x_1)^2 \,\le\, \frac{e^{\gamma T}}{\eta\gamma}\, 
  \widehat e_*(0)~, \qquad \hbox{for all } x_1 \in \R~,
\end{equation}
see \cite[Proposition~3.1]{GS2} for a similar argument. So, if 
we integrate \eqref{Fdiff1} over $x_1 \in [a{-}1,a{+}1]$ and 
then use \eqref{Fbound}, we arrive at the inequality
\[
  \,e^{\gamma T}\int_{a-1}^{a+1}\widehat e(x_1,T)\dd x_1 \,\le\, 
  2 \widehat e_*(0) + F(a{+}1) - F(a{-}1) \,\le\, 2 \widehat e_*(0)
  + 2\frac{e^{\gamma T/2}}{(\eta\gamma)^{1/2}}\,\widehat e_*(0)^{1/2}~,
\]
which gives \eqref{hatubd} since $\eta,\gamma$ are given by 
\eqref{bgdef} and $\widehat e_*(0) \le C\kappa^2$ with 
$\kappa < 1$. \QED

\bigskip
Proposition~\ref{hatudecay} shows that the oscillating part 
of the velocity $\widehat u(x,t)$ converges exponentially 
to zero as $t \to +\infty$ in the uniformly local space 
$L^2_\ul(\Omega)$, whose norm is defined as follows\:
\[
  \|\widehat u\|_{L^2_\ul(\Omega)}^2 \,=\, \sup_{a \in \R} 
  \int_{a-1}^{a+1} \!\!\int_\T |\widehat u(x_1,x_2)|^2 \dd x_2 
  \dd x_1~,
\]
see e.g. \cite{ARCD,GS1}. To conclude the proof of Theorem~\ref{main}, 
it remains to verify that we also have exponential decay in 
$L^\infty(\Omega)$. This follows directly from the following result\:

\begin{lemma}\label{smoothlem}
Assume that $u(x,t)$ is a solution of \eqref{NSeq2}, \eqref{pdef} 
satisfying
\begin{equation}\label{abound}
  \sup_{t \ge 0} \Bigl(\|u(t)\|_{L^\infty(\Omega)} + 
  \|\omega(t)\|_{L^\infty(\Omega)}\Bigr) \,\le\, M~, 
\end{equation}
for some $M > 0$. Then there exist positive constants $\tau, C_{10}$ 
such that, for all $t \ge 0$, 
\begin{equation}\label{smoothing}
  \|\widehat u(t+\tau)\|_{L^\infty(\Omega)} \,\le\, C_{10} 
  \|\widehat u(t)\|_{L^2_\ul(\Omega)}~.
\end{equation}
\end{lemma}

\noindent{\bf Proof.} We use equations \eqref{hatu1eq}, \eqref{hatu2eq}, 
which can be written in the compact form
\begin{equation}\label{wueq}
  \partial_t \widehat u + \partial_1 A_1(u) + \partial_2 A_2(u) 
  + B(u) \,=\, \Delta \widehat u~, 
\end{equation}
where
\[
  A_1 \,=\, \begin{pmatrix} \widehat u_1^2 + p \\ \widehat u_1 
  \widehat u_2 - \langle \widehat u_1\widehat u_2\rangle 
  \end{pmatrix}~, \qquad
  A_2 \,=\, \begin{pmatrix} (m+\widehat u_2)\widehat u_1 \\
  (m+\widehat u_2)\widehat u_2 +p \end{pmatrix}~, \qquad
  B = \begin{pmatrix} 0 \\ \widehat u_1 \partial_1 m\end{pmatrix}~.
\]
Here it is understood that the velocity field $u$ is 
decomposed as in \eqref{udecomp}, and that the pressure $p$ is 
given by \eqref{pdef}. The integral equation associated to 
\eqref{wueq} is
\begin{equation}\label{wuint}
  \widehat u(t) \,=\, S(t-t_0)\widehat u(t_0) - 
  \int_{t_0}^t \Bigl(\nabla \cdot S(t-s)A(u(s)) + S(t-s)B(u(s))\Bigr)\dd s~, 
\end{equation}
where $A = (A_1,A_2)^\top$ and $S(t) = e^{t\Delta}$ is the heat semigroup 
in $\Omega$. We have the following smoothing estimate
\begin{equation}\label{Sbdd}
  \|S(t)f\|_{L^\infty(\Omega)} \,\le\, \frac{C}{\min\{1,
  \sqrt{t}\}}\,\|f\|_{L^2_\ul(\Omega)}~, \qquad t > 0,
\end{equation}
which is easily established if we extend the function $f$ by 
periodicity and use the corresponding bound for the heat 
semigroup in the whole plane $\R^2$ \cite{ARCD}. On the 
other hand, in view of \eqref{abound} and \eqref{pdef}, 
we have the following bound on the nonlinear terms 
in \eqref{wuint}\:
\begin{equation}\label{ABbdd}
  \|A(u)\|_{L^\infty(\Omega)} + \|B(u)\|_{L^\infty(\Omega)} \,\le\, 
  C M \|\widehat u\|_{L^\infty(\Omega)}~.
\end{equation}

Now, we fix $t_0 \ge 0$ and assume that $t_0 < t \le t_0+1$. 
Using \eqref{Sbdd}, \eqref{ABbdd}, we obtain the following estimate 
for the solution of \eqref{wuint}\:
\[
  \|\widehat u(t)\|_{L^\infty} \,\le\, \frac{C}{(t{-}t_0)^{1/2}}
  \,\|\widehat u(t_0)\|_{L^2_\ul} + \int_{t_0}^t \frac{C'M}{(t{-}s)^{1/2}}
  \|\widehat u(s)\|_{L^\infty}\dd s~,
\]
for some $C,C' > 0$. In particular, if we denote $\Phi(t) = 
\sup\{(s{-}t_0)^{1/2}\|\widehat u(s)\|_{L^\infty}\,|\, t_0 < s \le t\}$, 
we see that
\begin{equation}\label{Phibdd}
  \Phi(t) \,\le\, C\|\widehat u(t_0)\|_{L^2_\ul} + \pi C'M 
  (t{-}t_0)^{1/2}\Phi(t)~, \qquad t_0 < t \le t_0 +1~.
\end{equation}
We now choose $\tau \in (0,1]$ such that $\pi C'M\tau^{1/2} \le 1/2$. 
It then follows from \eqref{Phibdd} that $\Phi(t) \le 
2C\|\widehat u(t_0)\|_{L^2_\ul}$ for $t_0 < t \le t_0+\tau$, hence
\[
  \|\widehat u(t_0+\tau)\|_{L^\infty} \,\le\, \frac{\Phi(t_0{+}
  \tau)}{\tau^{1/2}} \,\le\, \frac{2C}{\tau^{1/2}}\,
  \|\widehat u(t_0)\|_{L^2_\ul}~,
\]
which proves \eqref{smoothing}. \QED

\begin{corollary}\label{fincor}
Assume that $u(x,t)$ is a solution of \eqref{NSeq2}, \eqref{pdef} 
in $\Omega$ with bounded initial data satisfying \eqref{kappa}. 
Then 
\[
  u(x,t) \,=\, \begin{pmatrix} c \cr m(x_1,t) \end{pmatrix}
  + \widehat u(x,t)~, \qquad x \in \Omega~, \quad t \ge 0~,
\]
where $c \in \R$ is a constant, $m(x_1,t)$ is a solution of 
\eqref{meq}, and $\|\widehat u(t)\|_{L^\infty} = \OO(e^{-\gamma t})$ 
as $t \to \infty$ for any $\gamma < 2\pi^2$.  
\end{corollary}

Combining the results of Proposition~\ref{Unifom} and 
Corollary~\ref{fincor}, and returning to the original variables, 
we obtain \eqref{main3}. The proof of Theorem~\ref{main}
is now complete. \QED

\section{Conclusion and perspectives}\label{secConc}

In this final section, we briefly present some results obtained by 
S. Zelik \cite{Ze2} for the Navier-Stokes equations
in the whole plane $\R^2$, and we compare them to the conclusions
of Theorem~\ref{main} which hold when periodicity is assumed 
in one space direction. We first mention that, in \cite{Ze2}, the
following more general equation is considered\:
\[  
  \partial_t u + (u\cdot \nabla)u \,=\, \Delta u -\alpha u 
  -\nabla p + g~, \qquad \div u \,=\, 0~, 
\]
which includes an additional dissipation term $-\alpha u$ with
constant coefficient $\alpha \ge 0$, as well as a divergence-free
external force $g(x)$. However, in the spirit of the present work, we
only discuss here the results of \cite{Ze2} in the particular case
where $\alpha = 0$ and $g = 0$.

Let $L^2_\ul(\R^2)$ be the uniformly local $L^2$ space on $\R^2$ 
defined by the norm
\[
  \|f\|_{L^2_\ul} \,=\, \sup_{x \in\R^2} \biggl(\int_{|y-x| \le 1}
  |f(y)|^2 \dd y\biggr)^{1/2}~.
\]
If $u_0 \in L^2_\ul(\R^2)^2$ is divergence-free, it is known that 
the Navier-Stokes equations \eqref{NSeq2} have a unique global 
solution with initial data $u_0$, provided the pressure $p$ is 
given by the formula
\begin{equation}\label{Pdef}
  p \,=\, \sum_{i,j = 1}^2 R_i R_j (u_i u_j)~,
\end{equation}
where $R_1,R_2$ are the Riesz transforms on $\R^2$ \cite{GMS,MT}. 
This solution is smooth for positive times, and in particular 
the vorticity $\omega = \partial_1 u_2 - \partial_2 u_1$ is bounded
for all $t > 0$. Since we are mainly interested in the long-time 
behavior, we may thus assume without loss of generality that 
$\omega_0 = \curl u_0 \in L^\infty(\R^2)$.  

\begin{proposition}\label{Zelikprop} {\bf \cite{Ze2}}
Assume that $u_0 \in L^2_\ul(\R^2)^2$, $\div u_0 = 0$,  
and $\omega_0 = \curl u_0 \in L^\infty(\R^2)$. Then there exists
a constant $K \ge 1$ (depending only on $\|u_0\|_{L^2_\ul}$ and 
$\|\omega_0\|_{L^\infty}$) such that the solution of \eqref{NSeq2}, 
\eqref{Pdef} in $\R^2$ with initial data $u_0$ satisfies
\begin{equation}\label{ulR2}
  \sup_{x \in \R^2} \frac{1}{K t^2}\biggl(\int_{|y-x| \le K t^2}
  |u(y,t)|^2 \dd y\biggr)^{1/2} \,\le\, C \|u_0\|_{L^2_\ul}~,
\end{equation}
for all $t \ge 1$, where $C > 0$ is a universal constant. 
\end{proposition}

Remarkably enough, the proof of Proposition~\ref{Zelikprop} 
given in \cite[Section 7]{Ze2} does not use the viscous 
dissipation term $\Delta u$ in the Navier-Stokes equation. 
This means that estimate \eqref{ulR2} also holds for bounded
solutions of the Euler equations in $\R^2$, as long as these
solutions remain sufficiently smooth. In contrast, we emphasize
that the viscous dissipation was used in the proof of 
Theorem~\ref{main}, in particular in Section~\ref{subens}. 

As was observed in \cite{Ze2}, estimate \eqref{ulR2} is in some sense
optimal. For instance, if the initial velocity $u_0$ is constant and
nonzero, then $u(x,t) = u_0$ for all $x \in \R^2$ and all $t > 0$,
hence \eqref{ulR2} is sharp.  However one should observe that, in the
left-hand side of \eqref{ulR2}, averages are taken over very large
disks of radius $K t^2$, whereas in Sections~\ref{subener} and
\ref{subens} the corresponding domains (determined by the localization
function $\chi_\rho$) have a much smaller diameter, of order
$\sqrt{\beta t}$. The reason for this discrepancy is that, in the
cylinder $\Omega = \R \times \T$, it was easy to freeze the Galilean
invariance of the system and to assume, as in Section~\ref{subsecP2},
that the horizontal velocity has zero vertical average. As is shown in
Section~\ref{subfond}, this condition \eqref{zeromean} allows us to
prove that solutions of the vorticity equation \eqref{Veq2} behave
diffusively (in the horizontal direction) as $t \to \infty$, which
suggests that the diffusion length $\OO(\sqrt{t})$ is appropriate to
describe the spreading of solutions to the Navier-Stokes equation
\eqref{NSeq2} in that particular case. In the whole plane $\R^2$, the
situation is more complicated, and if we do not eliminate somehow
the Galilean invariance we are forced to take averages over disks of
radius at least $\OO(Ut)$, where $U$ is an upper bound on
$\|u\|_{L^\infty}$. In fact, since no a priori control on $U$ is 
available, the proof of Proposition~\ref{Zelikprop} is rather delicate
and relies on a ``self-consistent 
argument'' which eventually gives \eqref{ulR2}. As $Kt^2 \ge 1$ 
for $t \ge 1$, we immediately deduce from \eqref{ulR2} that 
$\|u(t)\|_{L^2_\ul} \le C K t^2 \|u_0\|_{L^2_\ul}$, but we do not 
know if that estimate is optimal. 

If we do use the viscous dissipation in the Navier-Stokes equations, 
then proceeding as in Section~\ref{secEner} it is possible to 
obtain the following result. 

\begin{corollary}\label{R2cor}
Under the assumptions of Proposition~\ref{Zelikprop}, the 
vorticity $\omega = \partial_1 u_2 - \partial_2 u_1$ satisfies, 
for all $t \ge 1$, 
\begin{equation}\label{ulomega}
  \sup_{x \in \R^2} \frac{1}{K t^2}\biggl(\int_{|y-x| \le K t^2}
  |\omega(y,t)|^2 \dd y\biggr)^{1/2} \,\le\, \frac{C}{\sqrt{t}} 
  \|u_0\|_{L^2_\ul}~.
\end{equation}
\end{corollary}

Estimate \eqref{ulomega} shows that the enstrophy of the solution, 
when averaged over sufficiently large disks, decays to zero as 
$t \to \infty$. This is the analog of Proposition~\ref{tElem}, 
except for the important discrepancy regarding the size of the disks,
which was already discussed. Unlike in the case of the cylinder, 
we are not able to convert \eqref{ulomega} into a uniform 
decay estimate for the vorticity. Nevertheless, estimate 
\eqref{ulomega} strongly suggests that the vorticity converges
to zero in some sense as $t \to \infty$, so that the long 
time asymptotics of \eqref{NSeq2} in $\R^2$ should be described
by irrotational flows, as was proved (in a particular case) 
in Theorem~\ref{main}. We hope to come back to this interesting
question in a future work. 

\section{Appendix}\label{secApp}

{\bf Proof of Lemma~\ref{intinlem}.} For any $p \in \{1\}\cup S$, 
it follows from \eqref{Fund2} that
\begin{equation}\label{Int1}
  w_{2p}'(t) \,\le\, -\frac{a}{p} \min_{\beta = 2,4}\biggl\{\Bigl(
  \frac{w_{2p}(t)}{w_p(t)} \Bigr)^{\beta p}\biggr\}w_{2p}(t) + pb\alpha^2
  w_{2p}(t)~, \qquad t > 0~,
\end{equation}
where $a = C/2$ and $b = 1+M^2/2$. We shall prove inductively that 
\begin{equation}\label{Int2}
  w_p(t) \,\le\, \overline w_p(t) \,:=\, \frac{A_p\,e^{B_p \alpha^2t}}{
  V(t)^{\frac{p-2}{2p}}}~, \qquad t > 0~, \quad p \in S~,
\end{equation}
provided the constants $A_p, B_p$ are chosen appropriately. First,
applying \eqref{Int1} with $p = 1$, we see that $w_2'(t) \,\le\,
\alpha^2 b\,w_2(t)$, hence \eqref{Int2} obviously holds for $p = 2$ if
$A_2 = w_2(0)$ and $B_2 = b$. Thus, it remains to show that, if
\eqref{Int2} holds for some $p \in S$, then the same inequality remains
true with $p$ replaced by $2p$, provided $A_{2p}$ and $B_{2p}$ are 
chosen appropriately.

To do that, we first observe that, if \eqref{Int2} holds for some 
$p \in S$, then the function $w_{2p}$ satisfies the differential
inequality
\begin{equation}\label{Int3}
  w_{2p}'(t) \,\le\, -\frac{a}{p} \min_{\beta = 2,4}\biggl\{\Bigl(
  \frac{w_{2p}(t)}{\overline w_p(t)} \Bigr)^{\beta p}\biggr\}w_{2p}(t) 
  + pb\alpha^2 w_{2p}(t)~, \qquad t > 0~,
\end{equation}
which is obtained from \eqref{Int1} by replacing $w_p(t)$ with
$\overline w_p(t)$. As we shall show, we can choose the constants 
$A_{2p}$, $B_{2p}$ so that the function $\overline w_{2p}$ defined 
by \eqref{Int2} satisfies the {\em reverse} inequality
\begin{equation}\label{Int4}
  \overline w_{2p}'(t) \,>\, -\frac{a}{p} \min_{\beta = 2,4}\biggl\{\Bigl(
  \frac{\overline w_{2p}(t)}{\overline w_p(t)} \Bigr)^{\beta p}\biggr\}
  \overline w_{2p}(t) + pb\alpha^2 \overline w_{2p}(t)~, \qquad t > 0~.
\end{equation}
Since obviously $\overline w_{2p}(t) \to +\infty$ as $t \to 0+$, it 
follows from \eqref{Int3}, \eqref{Int4} that $w_{2p}(t) \le 
\overline w_{2p}(t)$ for all $t > 0$, which proves \eqref{Int2}.  

It remains to establish \eqref{Int4}. Using the definition 
\eqref{Int2} of $\overline w_p$ and $\overline w_{2p}$, we find by a direct
calculation
\[
  \frac{\overline w_{2p}'(t)}{\overline w_{2p}(t)} \,=\, B_{2p}\alpha^2 - 
  \frac{p-1}{2p}\frac{V'(t)}{V(t)}~, \qquad \hbox{and}\quad
  \Bigl(\frac{\overline w_{2p}(t)}{\overline w_p(t)} \Bigr)^{\beta p}
  \,=\, \Bigl(\frac{A_{2p}}{A_p}\Bigr)^{\beta p} 
  \,\frac{e^{\beta p(B_{2p}-B_p)\alpha^2 t}}{V(t)^{\beta/2}}~.  
\]
Thus \eqref{Int4} holds provided
\begin{equation}\label{Int5}
  B_{2p}\alpha^2 - \frac{p-1}{2p}\frac{V'(t)}{V(t)} \,>\, 
  pb\alpha^2 - \frac{a}{p}\Bigl(\frac{A_{2p}}{A_p}\Bigr)^{\beta p} 
  \,\frac{e^{\beta p(B_{2p}-B_p)\alpha^2 t}}{V(t)^{\beta/2}}~, 
  \qquad t > 0~, 
\end{equation}
for $\beta = 2,4$. We now fix some $\epsilon \in (0,1)$ and choose 
$N \ge 1$ such that $2\epsilon a N \ge 1$. We assume that the 
constants $A_p, B_p$ in \eqref{Int2} satisfy the recursion relations
\begin{equation}\label{Int6}
  A_{2p} \,=\, A_p (Np^2)^{1/(2p)}~, \qquad \hbox{and}\quad
  B_{2p} \,=\, B_p(1 + \epsilon/p)~, \qquad p \in S~.
\end{equation}
Then \eqref{Int5} is equivalent to
\[
  B_{2p}\alpha^2 + \frac{a}{p}\Bigl(\frac{Np^2}{V(t)}\Bigr)^{\beta/2} 
  \,e^{\epsilon\beta B_p\alpha^2 t} \,>\, pb\alpha^2 
  + \frac{p-1}{2p}\frac{V'(t)}{V(t)}~, \qquad t > 0~. 
\]
In fact, since $B_p \ge B_2 = b$, we have $e^{\epsilon\beta B_p\alpha^2 t}
\ge 1 + \epsilon \beta b \alpha^2 t$, and it is thus sufficient to 
establish the stronger inequality
\[
  \frac{a}{p}\Bigl(\frac{Np^2}{V(t)}\Bigr)^{\beta/2} 
  \Bigl(1 + \epsilon \beta b \alpha^2 t\Bigr)\,\ge\, pb\alpha^2 
  + \frac{V'(t)}{2V(t)}~, \qquad t > 0~, 
\]
which is obviously satisfied if we can prove that
\begin{equation}\label{Int7}
  \frac{a}{p}\Bigl(\frac{Np^2}{V(t)}\Bigr)^{\beta/2} \,\ge\, 
  \frac{V'(t)}{2V(t)}~, \qquad \hbox{and} \qquad 
  \frac{a}{p}\Bigl(\frac{Np^2}{V(t)}\Bigr)^{\beta/2} \epsilon \beta t
  \,\ge\, p~, \qquad t > 0~. 
\end{equation}
But, since $V(t) = \min(t,\sqrt{t})$, it is clear that \eqref{Int7}
holds for $t > 0$ and $\beta = 2,4$ if $N \ge 1$ and $2\epsilon a N 
\ge 1$. This concludes the proof of the upper bound \eqref{Int2}. 

Finally, we iterate the recursion relations \eqref{Int6} to show
that the coefficients $A_p$, $B_p$ are uniformly bounded. A direct
calculation shows that
\[
  \sup_{p \in S} A_p \,=\, A_2 \prod_{p \in S}(p^2 N)^{1/(2p)} \,=\, 
  4A_2 N^{1/2}~, \qquad \hbox{and} \quad 
  \sup_{p \in S} B_p \,=\, B_2 \prod_{p \in S}(1+\epsilon/p) \,\le\, 
  B_2 e^\epsilon~.
\]
Since $A_2 = w_2(0)$ and $B_2 = b = 1+M^2/2 \le 1+M^2$, we see that 
\eqref{intineq} follows from \eqref{Int2}. \QED


\begin{thebibliography}{99}
\setlength{\itemsep}{-0.4mm}

\bibitem{AM} A. Afendikov and A. Mielke, Dynamical properties of
spatially non-decaying 2D Navier-Stokes flows with Kolmogorov 
forcing in an infinite strip, J. Math. Fluid. Mech. \textbf{7} 
(2005), suppl. 1, S51--S67.

\bibitem{AZ} P. Anthony and S. Zelik, 
Infinite-energy solutions for the Navier-Stokes equations in a strip 
revisited, preprint {\tt arXiv:1311.3128}. 

\bibitem{ARCD} J. Arrieta, A. Rodriguez-Bernal, J. Cholewa, and
T. Dlotko, Linear parabolic equations in locally uniform
spaces, Math. Models Methods Appl. Sci.  \textbf{14} (2004),
253--293.

\bibitem{Ar} D. G. Aronson, Bounds for the fundamental solution of 
a parabolic equation. Bull. Amer. Math. Soc. \textbf{73} (1967), 
890--896. 

\bibitem{BL}
J. Bergh and J. L\"ofstr\"om, {\em Interpolation spaces. An 
introduction.} Grundlehren der Mathematischen Wissenschaften 
\textbf{223}, Springer, 1976. 

\bibitem{CF} P. Constantin and C. Foias, 
{\em Navier-Stokes equations}, Chicago Lectures in Mathematics, 
University of Chicago Press, 1988. 

\bibitem{Co} Th. Coulhon, Ultracontractivity and Nash type inequalities.
J. Funct. Anal. \textbf{141} (1996), 510--539.

\bibitem{CGL} Th. Coulhon, A. Grigorʹyan, and D. Levin,
On isoperimetric profiles of product spaces.
Comm. Anal. Geom. \textbf{11} (2003), 85--120. 

\bibitem{FS} E. B. Fabes and D. W. Stroock, A new proof of Moser's 
parabolic Harnack inequality using the old ideas of Nash.
Arch. Rational Mech. Anal. \textbf{96} (1986), 327--338. 

\bibitem{GS1} Th. Gallay and S. Slijep\v{c}evi\'{c}, Energy flow in
formally gradient partial differential equations on unbounded domains, J.
Dynam. Differential Equations \textbf{13} (2001), 757--789.

\bibitem{GS2} Th. Gallay and S. Slijep\v{c}evi\'{c}, Distribution of
energy and convergence to equilibria in extended dissipative systems,
to appear in J.Dynam. Differential Equations. 

\bibitem{GS3} Th. Gallay and S. Slijep\v{c}evi\'{c}, Energy bounds 
for the two-dimensional Navier-Stokes equations in an infinite
cylinder, to appear in Commun. in PDE's. 

\bibitem{GMS} Y. Giga, S. Matsui, O. Sawada, Global existence of two
dimensional Navier-Stokes flow with non-decaying initial velocity, 
J. Math. Fluid. Mech. \textbf{3} (2001), 302--315.

\bibitem{Gr1} A. Grigorʹyan, 
Heat kernel upper bounds on a complete non-compact manifold,
Rev. Mat. Iberoamericana \textbf{10} (1994), 395--452. 

\bibitem{Gr2} A. Grigor'yan, 
{\em Heat kernel and analysis on manifolds}, 
AMS/IP Studies in Advanced Mathematics \textbf{47}, 
AMS, Providence, 2009.

\bibitem{Ka} J. Kato, The Uniqueness of Nondecaying Solutions 
for the Navier-Stokes Equations, Arch. Ration. Mech. Anal. 
\textbf{169} (2003), 159--175. 

\bibitem{MT} Y. Maekawa and Y. Terasawa,
The Navier-Stokes equations with initial data in uniformly local 
$L^p$ spaces, Diff. Int. Equations \text{bf 19} (2006), 369--400.

\bibitem{Na} J. Nash, Continuity of solutions of parabolic and 
elliptic equations, Amer. J. Math. \textbf{80} (1958), 931--954. 

\bibitem{ST} O. Sawada and Y. Taniuchi, A remark on $L^{\infty }$
solutions to the 2-D Navier-Stokes equations, J. Math. Fluid Mech. 
\textbf{9} (2007), 533--542.

\bibitem{Ze1} S. Zelik, Spatially nondecaying solutions of the 
2D Navier-Stokes equation in a strip, Glasg. Math. J. \textbf{49} 
(2007), 525--588. 

\bibitem{Ze2} S. Zelik, Infinite energy solutions for damped
Navier-Stokes equations in $\R^2$, J. Math. Fluid Mech. \textbf{15}
(2013), 717--745.

\end{thebibliography}
\end{document}